\documentclass[a4paper,twoside,10pt]{amsart}
\usepackage{amsmath,amsfonts,amssymb,amsthm}
\newtheorem{theorem}{Theorem}[section]

\usepackage{mathrsfs}
\usepackage{color,graphicx}
\usepackage[a4paper]{geometry}
\numberwithin{equation}{section}

\author[G. Nemes]{Gerg\H{o} Nemes}
\address{Central European University, Department of Mathematics and its Applications, H-1051 Budapest, N\'ador utca 9, Hungary}
\email{nemesgery@gmail.com}


\begin{document}

\title[Asymptotics of the Lommel function]{On the large argument asymptotics of\\ the Lommel function via Stieltjes transforms}

\begin{abstract} The aim of this paper is to investigate in detail the known large argument asymptotic series of the Lommel function by Stieltjes transform representations. We obtain a number of properties of this asymptotic expansion, including explicit and realistic error bounds, exponentially improved asymptotic expansions, and the smooth transition of the Stokes discontinuities. An interesting consequence related to the large argument asymptotic series of the Struve function is also proved.
\end{abstract}
\maketitle

\section{Introduction and main results}\label{section1}

In his important paper \cite{Boyd}, Boyd investigated various properties of the large-$z$ asymptotics of the modified Bessel function $K_\nu\left(z\right)$, using its Stieltjes transform representation. His analysis includes error bounds, exponentially improved asymptotic expansions and the smooth transition of the Stokes phenomenon.

In this paper, we shall discuss the Lommel function $S_{\mu ,\nu } \left( z \right)$ for fixed complex $\mu$, $\nu$ and large complex $z$. The asymptotic expansion for this function is well known \cite[11.9.iii]{NIST}. Following the notations in \cite{NIST}, we have that as $z \to \infty$ in the sector $\left|\arg z\right| \leq \pi -\delta$, for any $0< \delta \leq \pi$,
\begin{equation}\label{eq1}
S_{\mu ,\nu } \left( z \right) \sim z^{\mu  - 1} \sum\limits_{n = 0}^\infty  {\left( { - 1} \right)^n \frac{{a_n \left( { - \mu ,\nu } \right)}}{{z^{2n} }}},
\end{equation}
where
\[
a_n \left( {\mu ,\nu } \right) = \prod\limits_{k = 1}^n {\left( {\left( {\mu  + 2k - 1} \right)^2  - \nu ^2 } \right)}  = 2^{2n} \left( {\frac{{\mu  - \nu  + 1}}{2}} \right)_n \left( {\frac{{\mu  + \nu  + 1}}{2}} \right)_n ,
\]
with the Pochhammer symbol $\left( x \right)_n  = \Gamma \left( {x + n} \right)/\Gamma \left( x \right)$. If either of $\mu\pm\nu$ equals a positive odd integer, then the right-hand side of \eqref{eq1} terminates and represents $S_{\mu ,\nu } \left( z \right)$ exactly.

Before stating our results, we mention some applications to other special functions. The Struve functions are related to the Lommel function by
\begin{equation}\label{eq44}
\mathbf{K}_\nu  \left( z \right) = \frac{2^{1 - \nu }}{\sqrt \pi  \Gamma \left(\nu  + \frac{1}{2} \right)}S_{\nu ,\nu } \left( z \right)
\end{equation}
and
\begin{equation}\label{eq45}
\mathbf{M}_\nu  \left( z \right) =  \mp ie^{ \mp \frac{\pi }{2}i\nu } \mathbf{K}_\nu  \left( {ze^{ \pm \frac{\pi }{2}i} } \right) \mp \frac{2}{{\pi i}}e^{ \mp \pi i\nu } K_\nu  \left( z \right), \quad - \pi  \le  \pm \arg z \le \frac{\pi }{2}.
\end{equation}
The first formula is given in \cite[p. 444]{Dingle}, the second follows from \cite[11.2.E2, 11.2.E5 and 11.2.E6]{NIST}. The large argument asymptotic series of the Struve functions are well known (see, e.g., \cite[11.6.i]{NIST}). Explicit error bounds for the asymptotic expansion of $\mathbf{K}_\nu  \left( z \right)$ when $\nu$ is real and $z$ is positive, were obtained by Watson \cite[p. 333]{Watson}. The results of this paper provide error bounds for this expansion under more general circumstances. Dingle \cite[p. 445]{Dingle} gave exponentially improved versions of the asymptotic expansion of $\mathbf{K}_\nu  \left( z \right)$, although, his results were obtained by formal and interpretive, rather than rigorous, methods. Our analysis also provides a mathematically precise treatment of Dingle's formal expansions. Another family of special functions related to the Lommel function are the Anger--Weber-type functions:
\[
\mathbf{A}_\nu  \left( z \right) = \frac{{S_{0,\nu } \left( z \right) - \nu S_{ - 1,\nu } \left( z \right)}}{\pi },
\]
\[
\mathbf{J}_\nu  \left( z \right) - J_\nu  \left( z \right)= \frac{{\sin \left( {\pi \nu } \right)}}{\pi }\left( {S_{0,\nu } \left( z \right) - \nu S_{ - 1,\nu } \left( z \right)} \right)
\]
and
\[
\mathbf{E}_\nu  \left( z \right) + Y_\nu  \left( z \right)= \frac{{\left( {\cos \left( {\pi \nu } \right) - 1} \right)\nu S_{ - 1,\nu } \left( z \right) - \left( {\cos \left( {\pi \nu } \right) + 1} \right)S_{0,\nu } \left( z \right)}}{\pi } .
\]
For the definitions and the large argument asymptotics, see, e.g., \cite[11.10.i and 11.11.i]{NIST}. These connection formulas follow by applying formula \cite[11.9.E5]{NIST} to \cite[11.10.E17]{NIST} and \cite[11.10.E18]{NIST} (see also Luke \cite[p. 84]{Luke}). Precise error bounds for the large argument asymptotic series of these functions were derived by Meijer \cite{Meijer}. The error bounds we prove in this paper are generalisations of Meijer's results.

In our first theorem, we give two Stieltjes transform-type integral representations for the remainder of the asymptotic series \eqref{eq1}. Throughout this paper, empty sums are taken to be zero.

\begin{theorem}\label{thm1} For complex $z$, $\mu$ and $\nu$, and for any non-negative integer $N$, define the remainder term $R_N \left( {z,\mu ,\nu } \right)$ by
\begin{equation}\label{eq27}
S_{\mu ,\nu } \left( z \right) = z^{\mu  - 1} \sum\limits_{n = 0}^{N - 1} {\left( { - 1} \right)^n \frac{a_n \left( { - \mu ,\nu } \right)}{z^{2n}}}  + R_N \left( {z,\mu ,\nu } \right) .
\end{equation}
Then we have the integral representations
\begin{equation}\label{eq2}
R_N \left( {z,\mu ,\nu } \right) = \left( { - 1} \right)^N \frac{{2^{\mu  + 1} z^{\mu  - 2N - 1} }}{{\Gamma \left( {\frac{{\nu  - \mu  + 1}}{2}} \right)\Gamma \left( {\frac{{1 - \mu  - \nu }}{2}} \right)}}\int_0^{ + \infty } {\frac{{t^{2N - \mu } K_\nu  \left( t \right)}}{{1 + \left( {t/z} \right)^2 }}dt} 
\end{equation}
for $\left|\arg z\right| < \frac{\pi}{2}$ and $\Re \left( \mu  \right) + \left| {\Re \left( \nu  \right)} \right|< 2N + 1$;
\begin{equation}\label{eq5}
R_N \left( {z,\mu ,\nu } \right) = 2^\mu  \frac{{\Gamma \left( {\frac{{\mu  + \nu  + 1}}{2}} \right)}}{{\Gamma \left( {\frac{{\nu  - \mu  + 1}}{2}} \right)}}\left( {z^{\mu  - 2N - 1} \int_0^{ + \infty } {\frac{{t^{2N - \mu } J_\nu  \left( t \right)}}{{1 - \left( {t/z} \right)^2 }}dt}  + e^{ \mp \frac{\pi }{2}i\nu } K_\nu  \left( {ze^{ \mp \frac{\pi }{2}i} } \right)} \right)
\end{equation}
for $0 < \pm \arg z <\pi$, $\Re \left( \mu  \right) - \Re \left( \nu  \right)< 2N + 1$ and $2N - \frac{3}{2} < \Re \left( \mu  \right)$.
\end{theorem}

In our analysis of the asymptotic expansion \eqref{eq1}, we shall use the representation \eqref{eq2}. Nevertheless, the formula \eqref{eq5} has some important consequences on the asymptotics of the Struve function $\mathbf{M}_\nu  \left( z \right)$ (see Section \ref{section5}). In Section \ref{section3}, we will show how to obtain numerically computable bounds for the remainder term $R_N \left( {z,\mu ,\nu } \right)$.

In his paper \cite{Boyd}, Boyd derived the representation
\begin{equation}\label{eq28}
K_\nu  \left( z \right) = \left( {\frac{\pi }{2z}} \right)^{\frac{1}{2}} e^{ - z} \sum\limits_{n = 0}^{N - 1} {\frac{{a_n\left( \nu  \right)}}{{z^n }}}  + K_N \left( {z,\nu } \right)
\end{equation}
for $\left|\arg z\right| < \pi$, with
\[
a_n \left( \nu  \right) = \left( { - 1} \right)^n \frac{{\left( {\frac{1}{2} + \nu } \right)_n \left( {\frac{1}{2} - \nu } \right)_n }}{{2^n \Gamma \left( {n + 1} \right)}} = \left( { - 1} \right)^n \sqrt {\frac{2}{\pi }} \frac{{\cos \left( {\nu \pi } \right)}}{\pi }\int_0^{ + \infty } {t^{n - \frac{1}{2}}e^{ - t} K_\nu  \left( t \right) dt} .
\]
The error term $K_N\left( {z,\nu } \right)$ is given by
\begin{equation}\label{eq29}
K_N \left( {z,\nu } \right) = \left( { - 1} \right)^N \frac{{\cos \left( {\pi \nu} \right)}}{\pi }z^{ - N - \frac{1}{2}} e^{ - z} \int_0^{ + \infty } {\frac{{t^{N - \frac{1}{2}} e^{ - t} K_\nu  \left( t \right)}}{{1 + t/z}}dt} ,
\end{equation}
provided that $\left| {\Re \left( \nu  \right)} \right| < N + \frac{1}{2}$. Employing \eqref{eq27}, \eqref{eq5}, \eqref{eq28}, \eqref{eq29}, the connection formula between the Lommel function and the Bessel functions (see Watson \cite[p. 347, expression (3)]{Watson} or \cite[11.9.E5]{NIST}), together with the continuation formulas for $K_\nu  \left( z \right)$ (see \cite[10.34.E4]{NIST}), it is possible to derive expansions for $S_{\mu ,\nu } \left( z \right)$ in other sectors of the complex $z$-plane.

In the following theorem, we give exponentially improved asymptotic expansion for the function $S_{\mu ,\nu } \left( z \right)$. This expansion can be viewed as the mathematically rigorous form of the terminated expansions of Dingle \cite[pp. 442--444]{Dingle}. In this theorem, we truncate the asymptotic series of $S_{\mu ,\nu } \left( z \right)$ at about its least term and re-expand the remainder into a new asymptotic expansion. The resulting exponentially improved asymptotic series is valid in larger regions than the original expansion \eqref{eq1}. The terms in this new series involve the Terminant function $\widehat T_p\left(w\right)$, which allows the smooth transition through the Stokes lines $\arg \nu = \pm \frac{\pi}{2}$. For the definition and basic properties of the Terminant function, see Section \ref{section4}. Throughout this paper, we use subscripts in the $\mathcal{O}$ notations to indicate the dependence of the implied constant on certain parameters.

\begin{theorem}\label{thm2} Let $M$ be an arbitrary fixed non-negative integer, and let $\mu$ and $\nu$ be fixed complex numbers. Suppose that $\left|\arg z\right| \leq \frac{3\pi}{2}$, $\left|z\right|$ is large and $N = \frac{1}{2}\left|z\right|+\rho$ is a positive integer with $\rho$ being bounded. Then
\begin{multline}\label{eq21}
R_N \left( {z,\mu ,\nu } \right) = \frac{{2^{\mu  + 1} }}{{\Gamma \left( {\frac{{\nu  - \mu  + 1}}{2}} \right)\Gamma \left( {\frac{{1 - \mu  - \nu }}{2}} \right)}}\left( {\pi e^{\frac{\pi }{2}i\mu } \left( {\frac{\pi}{-2iz}} \right)^{\frac{1}{2}} e^{iz} \sum\limits_{m = 0}^{M - 1} {\frac{{a_m \left( \nu  \right)}}{{\left( {-iz} \right)^m }} \widehat T_{2N - m - \mu  + \frac{1}{2}} \left( {iz} \right)} }\right. \\ \left.{+ \pi e^{ - \frac{\pi }{2}i\mu } \left( {\frac{\pi }{{2iz }}} \right)^{\frac{1}{2}} e^{ - iz} \sum\limits_{m = 0}^{M - 1} {\frac{{a_m \left( \nu  \right)}}{{\left( {iz } \right)^m }}e^{2\pi i\mu } \widehat T_{2N - m - \mu  + \frac{1}{2}} \left( { - iz} \right)}  + R_{N,M} \left( {z,\mu ,\nu } \right)} \right),
\end{multline}
where
\[
R_{N,M} \left( {z,\mu ,\nu } \right) = \mathcal{O}_{M,\mu ,\rho } \left( {\frac{{e^{ - \left| z \right|} }}{{\left| z \right|^{\frac{1}{2}} }}\frac{{\left| {a_M \left( \nu  \right)} \right|}}{{\left| z \right|^M }}} \right)
\]
for $\left|\arg z\right| \leq \frac{\pi}{2}$;
\[
R_{N,M} \left( {z,\mu ,\nu } \right) = \mathcal{O}_{M,\mu ,\rho } \left( {\frac{{e^{ \mp \Im \left( z \right)} }}{{\left| z \right|^{\frac{1}{2}} }}\frac{{\left| {a_M \left( \nu  \right)} \right|}}{{\left| z \right|^M }}} \right)
\]
for $\frac{\pi }{2} \leq \pm \arg z \leq \frac{{3\pi }}{2}$.
\end{theorem}

While proving Theorem \ref{thm2} in Section \ref{section4}, we also obtain the following explicit bound for the remainder $R_{N,M} \left( {z,\mu ,\nu } \right)$ in \eqref{eq21}. Note that in this theorem $N$ may not depend on $\left|z\right|$.

\begin{theorem}\label{thm3} Let $N$ and $M$ be non-negative integers and let $\mu$ and $\nu$ be complex numbers, such that $\left| {\Re \left( \nu  \right)} \right| < M + \frac{1}{2}$ and $\Re \left( \mu  \right) < 2N - M + \frac{1}{2}$. Define the remainder $R_{N,M} \left( {z,\mu ,\nu } \right)$ by the equality \eqref{eq21}. Then we have
\begin{align*}
\left| {R_{N,M} \left( {z,\mu ,\nu } \right)} \right| \le \; & \pi \left| {e^{\frac{\pi }{2}i\mu } } \right|\left( {\frac{\pi }{{2\left| z \right|}}} \right)^{\frac{1}{2}} \frac{{\left| {\cos \left( {\pi \nu } \right)} \right|}}{{\left| {\cos \left( {\pi \Re \left( \nu  \right)} \right)} \right|}}\frac{{\left| {a_M \left( \Re \left( \nu  \right) \right)} \right|}}{{\left| z \right|^M }}\left| {e^{iz} \widehat T_{2N - M - \mu  + \frac{1}{2}} \left( {iz} \right)} \right| \\ & + \pi \left| {e^{ - \frac{\pi }{2}i\mu } } \right|\left( {\frac{\pi }{{2\left| z \right|}}} \right)^{\frac{1}{2}} \frac{{\left| {\cos \left( {\pi \nu } \right)} \right|}}{{\left| {\cos \left( {\pi \Re \left( \nu  \right)} \right)} \right|}}\frac{{\left| {a_M \left( \Re \left( \nu  \right) \right)} \right|}}{{\left| z \right|^M }}\left| {e^{2\pi i\mu } e^{ - iz} \widehat T_{2N - M - \mu  + \frac{1}{2}} \left( { - iz} \right)} \right|
\\ & + \left( {\frac{\pi }{2}} \right)^{\frac{1}{2}} \left| {z^\mu  } \right|\frac{{\left| {\cos \left( {\pi \nu } \right)} \right|}}{{\left| {\cos \left( {\pi \Re \left( \nu  \right)} \right)} \right|}}\frac{{\left| {a_M \left( \Re \left( \nu  \right) \right)} \right|\Gamma \left( {2N - M - \Re \left( \mu  \right) + \frac{1}{2}} \right)}}{{\left| z \right|^{2N + 1} }},
\end{align*}
provided that $\left|\arg z\right| \leq \frac{\pi}{2}$. When $2\Re\left(\nu\right)$ is an odd integer, the limiting value has to be taken in this bound.
\end{theorem}

In their paper \cite{Howls}, Howls and Olde Daalhuis investigated the hyperasymptotic properties of solutions of inhomogeneous linear differential equations with a singularity of rank one. The result in Theorem \ref{thm2} can be regarded as a special case of their theory. Nevertheless, our approach provides not only an order estimate but an explicit, numerically computable bound for the remainder $R_{N,M} \left( {z,\mu ,\nu } \right)$.

If we restrict $z$ to the right-half plane, the re-expansion of the remainder $R_{N} \left( {z,\mu ,\nu } \right)$ can be done using only elementary functions. For a general theory of such re-expansions, see the papers of Olde Daalhuis \cite{Olde Daalhuis1,Olde Daalhuis2}.

\begin{theorem}\label{thm4} Let $M$ be an arbitrary fixed non-negative integer. Let $\mu$ and $\nu$ be fixed complex numbers. Suppose that $\left|\arg z\right| \leq \frac{\pi}{2}-\delta < \frac{\pi}{2}$ with $0<\delta \leq \frac{\pi}{2}$ being fixed, $\left|z\right|$ is large and $N = \frac{1}{2}\left|z\right|+\rho$ is a positive integer with $\rho$ being bounded. Then
\begin{equation}\label{eq42}
R_{N} \left( {z,\mu ,\nu } \right) = z^{\mu - 1} \sum\limits_{m = 0}^{M - 1} {v_{N,m} \left( {z,\mu ,\nu } \right)}  + \widehat R_{N,M} \left( {z,\mu ,\nu } \right),
\end{equation}
where
\begin{align*}
v_{N,m} \left( {z,\mu ,\nu } \right) & = \frac{\left( { - 1} \right)^N}{z^{2N}} \frac{{e^{ - 2i\theta m} }}{{\left( {1 + e^{ - 2i\theta } } \right)^{m + 1} }}\sum\limits_{k = 0}^m {\binom{m}{k}\left( { - 1} \right)^k \frac{{a_{N + k} \left( { - \mu ,\nu } \right)}}{{\left| z \right|^{2k} }}} 
\\ & = \frac{1}{{\left( {1 + e^{ - 2i\theta } } \right)^{m + 1} }}\sum\limits_{k = 0}^m {\binom{m}{k}\left( { - 1} \right)^{N + k} \frac{{a_{N + k} \left( { - \mu ,\nu } \right)}}{{z^{2N + 2k} }}e^{ - 2i\theta \left( {m - k} \right)} }
\end{align*}
and
\begin{equation}\label{eq52}
\widehat R_{N,M} \left( {z,\mu ,\nu } \right) = \mathcal{O}_{M,\mu ,\nu ,\rho,\delta } \left( {\frac{{e^{ - \left| z \right|} }}{{\left| z \right|^{\frac{M}{2} + 1} }}} \right).
\end{equation}
\end{theorem}

Note that the expansion \eqref{eq42} may be regarded as the result of an Euler transformation \cite[p. 537]{Olver} applied to the divergent tail of the asymptotic series \eqref{eq1}.

For some explicit bounds for the remainder $\widehat R_{N,M} \left( {z,\mu ,\nu } \right)$ in \eqref{eq42}, when $M$ is even and $N$ may not depend on $\left|z\right|$, see Subsection \ref{subsection43}.

The idea of re-expanding the optimally truncated remainder into another asymptotic series, in order to improve its numerical efficacy, dates back to Stieltjes \cite{Stieltjes}. In Stieltjes' work, the object of re-expansion is not the remainder but the converging factor, which is the ratio of the remainder and the first omitted term. For the sake of interest, we compute the asymptotic series of the converging factor corresponding to the asymptotic expansion \eqref{eq1} of the Lommel function in Appendix \ref{appendixa}. For more information about converging factors, see Olver \cite{Olver3}, \cite[pp. 522--536]{Olver} or Paris \cite[pp. 85--89]{Paris}.

The rest of the paper is organised as follows. In Section \ref{section2}, we prove the formulas for the remainder term stated in Theorem \ref{thm1}. In Section \ref{section3}, we give explicit and numerically computable error bounds for the large argument asymptotic series \eqref{eq1} of the Lommel function using the results in Theorem \ref{thm1}. In Section \ref{section4}, we prove the exponentially improved expansions presented in Theorems \ref{thm2} and \ref{thm4} together with the error bound in Theorem \ref{thm3}, and provide a detailed discussion of the Stokes phenomenon related to the expansion \eqref{eq1}. An interesting property of the large argument asymptotic series of the Struve function $\mathbf{M}_\nu  \left( z \right)$, which was already observed by Dingle, is proved in Section \ref{section5}.

\section{Proof of the formulas for the remainder term}\label{section2}

First, we prove \eqref{eq2}. Our starting point is the integral representation (see Dingle \cite[p. 380]{Dingle} or Erd\'elyi et al. \cite[p. 230, entry (42)]{Erdelyi2})
\begin{equation}\label{eq3}
S_{\mu ,\nu } \left( z \right) = \frac{{2^{\mu  + 1} z^{\mu  - 1} }}{{\Gamma \left( {\frac{{\nu  - \mu  + 1}}{2}} \right)\Gamma \left( {\frac{{1 - \mu  - \nu }}{2}} \right)}}\int_0^{ + \infty } {\frac{{t^{ - \mu } K_\nu  \left( t \right)}}{{1 + \left( {t/z} \right)^2 }}dt} ,
\end{equation}
valid for $\left|\arg z\right| < \frac{\pi}{2}$ and $\Re \left( \mu  \right) + \left| {\Re \left( \nu  \right)} \right|<1$. For every non-negative integer $N$ and $\left|\arg z\right| < \frac{\pi}{2}$ we have
\[
\frac{1}{{1 + \left( {t/z} \right)^2 }} = \sum\limits_{n = 0}^{N - 1} {\frac{{\left( { - 1} \right)^n }}{{z^{2n} }}t^{2n} }  + \frac{{\left( { - 1} \right)^N }}{z^{2N}}\frac{{t^{2N} }}{{1 + \left( {t/z} \right)^2 }} .
\]
Substituting this into the expression \eqref{eq3} and integrating term-by-term, we deduce
\begin{gather}\label{eq4}
\begin{split}
S_{\mu ,\nu } \left( z \right) = \; & z^{\mu  - 1} \sum\limits_{n = 0}^{N - 1} {\frac{{\left( { - 1} \right)^n }}{{z^{2n} }}\frac{{2^{\mu  + 1} }}{{\Gamma \left( {\frac{{\nu  - \mu  + 1}}{2}} \right)\Gamma \left( {\frac{{1 - \mu  - \nu }}{2}} \right)}}\int_0^{ + \infty } {t^{2n - \mu } K_\nu  \left( t \right)dt} } \\ & + \left( { - 1} \right)^N \frac{{2^{\mu  + 1} z^{\mu  - 2N - 1} }}{{\Gamma \left( {\frac{{\nu  - \mu  + 1}}{2}} \right)\Gamma \left( {\frac{{1 - \mu  - \nu }}{2}} \right)}}\int_0^{ + \infty } {\frac{{t^{2N - \mu } K_\nu  \left( t \right)}}{{1 + \left( {t/z} \right)^2 }}dt} .
\end{split}
\end{gather}
As $t\to 0^+$, we have $K_\nu  \left( t \right) = \mathcal{O}\left( {t^{ - \left|\Re\left(\nu\right)\right| } } \right)$ for $\nu \neq 0$, and $K_\nu  \left( t \right) = \mathcal{O}\left(\log t \right)$ for $\nu = 0$. Also $K_\nu  \left( t \right) = \mathcal{O}\left( {t^{ - \frac{1}{2}} e^{ - t} } \right)$ as $t\to +\infty$. Therefore, the integrals in \eqref{eq4} converge as long as $\Re \left( \mu  \right) + \left| {\Re \left( \nu  \right)} \right|< 2N + 1$; and by analytic continuation, \eqref{eq4} is valid for $\left|\arg z\right| < \frac{\pi}{2}$ and $\Re \left( \mu  \right) + \left| {\Re \left( \nu  \right)} \right|< 2N + 1$. The integrals under the sum evaluate to
\[
\frac{{2^{\mu  + 1} }}{{\Gamma \left( {\frac{{\nu  - \mu  + 1}}{2}} \right)\Gamma \left( {\frac{{1 - \mu  - \nu }}{2}} \right)}}\int_0^{ + \infty } {t^{2n - \mu } K_\nu  \left( t \right)dt}  = 2^{2n} \frac{{\Gamma \left( {\frac{{\nu  - \mu  + 1}}{2} + n} \right)\Gamma \left( {\frac{{1 - \mu  - \nu }}{2} + n} \right)}}{{\Gamma \left( {\frac{{\nu  - \mu  + 1}}{2}} \right)\Gamma \left( {\frac{{1 - \mu  - \nu }}{2}} \right)}} = a_n \left( { - \mu ,\nu } \right)
\]
(see Erd\'elyi et al. \cite[p. 331, entry (26)]{Erdelyi1}), which completes the proof of \eqref{eq2}.

Next we prove \eqref{eq5}. We start with the formula (see Erd\'elyi et al. \cite[p. 225, entry (13)]{Erdelyi2})
\begin{multline*}
{}_1F_2 \left( {1;1 - \lambda  - \frac{\nu }{2},1 - \lambda  + \frac{\nu }{2};\frac{z^2}{4}} \right) - \Gamma \left( {1 - \lambda  - \frac{\nu }{2}} \right)\Gamma \left( {1 - \lambda  + \frac{\nu }{2}} \right)\left( {\frac{z}{2}} \right)^{2\lambda } I_\nu  \left( z \right) \\ = \frac{{\pi 2^{1 - 2\lambda }  z^{ - 2} }}{{\sin \left( {\left( {\lambda  - \frac{\nu }{2}} \right)\pi } \right)\Gamma \left( {\lambda  + \frac{\nu }{2}} \right)\Gamma \left( {\lambda  - \frac{\nu }{2}} \right)}}\int_0^{ + \infty } {\frac{{t^{2\lambda  + 1} J_\nu  \left( t \right)}}{{1 + \left( {t/z} \right)^2 }}dt} ,
\end{multline*}
valid for $\Re \left( {\lambda  + \frac{\nu }{2}} \right) >  - 1$, $\Re \left( \lambda  \right) < \frac{1}{4}$ and $\left|\arg z\right| <\frac{\pi}{2}$. Substituting $\lambda  =  - \frac{\mu }{2} - \frac{1}{2}$ and $z = ze^{ \pm \frac{\pi }{2}i}$, we arrive at
\begin{multline*}
{}_1F_2 \left( {1;\frac{{\mu  - \nu  + 3}}{2},\frac{{\mu  + \nu  + 3}}{2}; - \frac{z^2}{4}} \right) - \Gamma \left( {\frac{{\mu  - \nu  + 3}}{2}} \right)\Gamma \left( {\frac{{\mu  + \nu  + 3}}{2}} \right)\left( {\frac{{ze^{ \mp \frac{\pi }{2}i} }}{2}} \right)^{ - \mu  - 1} I_\nu  \left( {ze^{ \mp \frac{\pi }{2}i} } \right) \\ =  \frac{{\pi 2^{2 + \mu } z^{ - 2} }}{{\cos \left( {\frac{{\mu  + \nu }}{2}\pi } \right)\Gamma \left( {\frac{{\nu  - \mu  - 1}}{2}} \right)\Gamma \left( {\frac{{ - \mu  - \nu  - 1}}{2}} \right)}}\int_0^{ + \infty } {\frac{{t^{ - \mu } J_\nu  \left( t \right)}}{{1 - \left( {t/z} \right)^2 }}dt} ,
\end{multline*}
provided that $\Re \left( \nu  \right) - \Re \left( \mu  \right) >  - 1$, $\Re \left( \mu  \right) >  - \frac{3}{2}$, $0 <  \pm \arg z < \pi$. This expression can be simplified using the Lommel function $s_{\mu ,\nu } \left( z \right)$ (see Watson \cite[p. 346, expression (10)]{Watson}) and the connection formula $I_\nu  \left( {ze^{ \mp \frac{\pi }{2}i} } \right) = e^{ \mp \frac{\pi}{2} i \nu } J_\nu  \left( z \right)$ for $0 <  \pm \arg z < \pi$, to obtain
\begin{multline*}
s_{\mu ,\nu } \left( z \right) - 2^{\mu  - 1} \Gamma \left( {\frac{{\mu  - \nu  + 1}}{2}} \right)\Gamma \left( {\frac{{\mu  + \nu  + 1}}{2}} \right)e^{ \pm \frac{\pi }{2}i\left( {\mu  - \nu  + 1} \right)} J_\nu  \left( z \right)
\\ = \frac{{\pi 2^\mu  z^{\mu  - 1} }}{{\cos \left( {\frac{{\mu  + \nu }}{2}\pi } \right)\Gamma \left( {\frac{{\nu  - \mu  + 1}}{2}} \right)\Gamma \left( {\frac{{1 - \mu  - \nu }}{2}} \right)}}\int_0^{ + \infty } {\frac{{t^{ - \mu } J_\nu  \left( t \right)}}{{1 - \left( {t/z} \right)^2 }}dt} .
\end{multline*}
Employing the expression
\begin{multline*}
 - e^{ \pm \frac{\pi }{2}i\left( {\mu  - \nu  + 1} \right)} J_\nu  \left( z \right) \\ = \sin \left( {\frac{{\mu  - \nu }}{2}\pi } \right)J_\nu  \left( z \right) - \cos \left( {\frac{{\mu  - \nu }}{2}\pi } \right)Y_\nu  \left( z \right) - \frac{2}{\pi }\cos \left( {\frac{{\mu  - \nu }}{2}\pi } \right)e^{ \mp \frac{\pi }{2}i\nu } K_\nu  \left( {ze^{ \mp \frac{\pi }{2}i} } \right)
\end{multline*}
and the connection formula between the functions $s_{\mu ,\nu } \left( z \right)$ and $S_{\mu ,\nu } \left( z \right)$ (see Watson \cite[p. 347, expression (3)]{Watson} or \cite[11.9.E5]{NIST}), we deduce
\begin{multline*}
S_{\mu ,\nu } \left( z \right) - \frac{{2^\mu  }}{\pi }\Gamma \left( {\frac{{\mu  - \nu  + 1}}{2}} \right)\Gamma \left( {\frac{{\mu  + \nu  + 1}}{2}} \right)\cos \left( {\frac{{\mu  - \nu }}{2}\pi } \right)e^{ \mp \frac{\pi }{2}i\nu } K_\nu  \left( {ze^{ \mp \frac{\pi }{2}i} } \right) \\ = \frac{{\pi 2^\mu  z^{\mu  - 1} }}{{\cos \left( {\frac{{\mu  + \nu }}{2}\pi } \right)\Gamma \left( {\frac{{\nu  - \mu  + 1}}{2}} \right)\Gamma \left( {\frac{{1 - \mu  - \nu }}{2}} \right)}}\int_0^{ + \infty } {\frac{{t^{ - \mu } J_\nu  \left( t \right)}}{{1 - \left( {t/z} \right)^2 }}dt} .
\end{multline*}
This can be simplified to the more compact form
\begin{equation}\label{eq6}
S_{\mu ,\nu } \left( z \right) = 2^\mu  \frac{{\Gamma \left( {\frac{{\mu  + \nu  + 1}}{2}} \right)}}{{\Gamma \left( {\frac{{\nu  - \mu  + 1}}{2}} \right)}}\left( {z^{\mu  - 1} \int_0^{ + \infty } {\frac{{t^{ - \mu } J_\nu  \left( t \right)}}{{1 - \left( {t/z} \right)^2 }}dt}  + e^{ \mp \frac{\pi }{2}i\nu } K_\nu  \left( {ze^{ \mp \frac{\pi }{2}i} } \right)} \right)
\end{equation}
valid for $\Re \left( \nu  \right) - \Re \left( \mu  \right) >  - 1$, $\Re \left( \mu  \right) >  - \frac{3}{2}$ and $0 <  \pm \arg z < \pi$. For every non-negative integer $N$ and $0 <  \pm \arg z < \pi$ we have
\[
\frac{1}{{1 - \left( {t/z} \right)^2 }} = \sum\limits_{n = 0}^{N - 1} {\frac{{t^{2n} }}{{z^{2n} }}}  + \frac{1}{{z^{2N} }}\frac{{t^{2N} }}{{1 - \left( {t/z} \right)^2 }}.
\]
Substituting this into the expression \eqref{eq6} and integrating term-by-term, we deduce
\begin{gather}\label{eq7}
\begin{split}
S_{\mu ,\nu } \left( z \right) = \; & z^{\mu  - 1} \sum\limits_{n = 0}^{N - 1} {\frac{1}{{z^{2n} }}2^\mu  \frac{{\Gamma \left( {\frac{{\mu  + \nu  + 1}}{2}} \right)}}{{\Gamma \left( {\frac{{\nu  - \mu  + 1}}{2}} \right)}}\int_0^{ + \infty } {t^{2n - \mu } J_\nu  \left( t \right)dt} } 
\\ & + 2^\mu  \frac{{\Gamma \left( {\frac{{\mu  + \nu  + 1}}{2}} \right)}}{{\Gamma \left( {\frac{{\nu  - \mu  + 1}}{2}} \right)}}\left( {z^{\mu  - 2N - 1} \int_0^{ + \infty } {\frac{{t^{2N - \mu } J_\nu  \left( t \right)}}{{1 - \left( {t/z} \right)^2 }}dt}  + e^{ \mp \frac{\pi }{2}i\nu } K_\nu  \left( {ze^{ \mp \frac{\pi }{2}i} } \right)} \right) .
\end{split}
\end{gather}
As $t\to 0^+$, we have $J_\nu  \left( t \right) = \mathcal{O}\left( t^\nu \right)$ for $\nu \neq -1,-2,-3,\ldots$, and $J_\nu  \left( t \right) = \mathcal{O}\left( t^{-\nu} \right)$ for $\nu = -1,-2,-3,\ldots$. Also $J_\nu  \left( t \right) = \mathcal{O}\left( t^{ - \frac{1}{2} } \right)$ as $t\to +\infty$. Therefore, the integrals in \eqref{eq7} converge as long as $\Re \left( \mu  \right) - \Re \left( \nu  \right)< 2N + 1$ and $2N - \frac{3}{2} < \Re \left( \mu  \right)$; and by analytic continuation, \eqref{eq7} is valid for $0 < \pm \arg z <\pi$, $\Re \left( \mu  \right) - \Re \left( \nu  \right)< 2N + 1$ and $2N - \frac{3}{2} < \Re \left( \mu  \right)$. The integrals under the sum evaluate to
\[
2^\mu  \frac{{\Gamma \left( {\frac{{\mu  + \nu  + 1}}{2}} \right)}}{{\Gamma \left( {\frac{{\nu  - \mu  + 1}}{2}} \right)}}\int_0^{ + \infty } {t^{2n - \mu } J_\nu  \left( t \right)dt}  = 2^{2n} \frac{{\Gamma \left( {\frac{{\nu  - \mu  + 1}}{2} + n} \right)\Gamma \left( {\frac{{\mu  + \nu  + 1}}{2}} \right)}}{{\Gamma \left( {\frac{{\nu  - \mu  + 1}}{2}} \right)\Gamma \left( {\frac{{\mu  + \nu  + 1}}{2} - n} \right)}} = \left( { - 1} \right)^n a_n \left( { - \mu ,\nu } \right)
\]
(see Erd\'elyi et al. \cite[p. 326, entry (1)]{Erdelyi1}), which completes the proof of \eqref{eq5}.

\section{Error bounds for the large argument asymptotics of the Lommel function}\label{section3}

In this section, we derive numerically computable bounds for the error term $R_N \left( {z,\mu ,\nu } \right)$ in \eqref{eq27}. We may assume that neither of $\mu\pm\nu$ equals an positive odd integer, otherwise the asymptotic series terminates and represents $S_{\mu ,\nu } \left( z \right)$ exactly. First we consider bounds which are suitable for $\Re\left(z\right)>0$ and when $z$ is not too close to the imaginary axis. To make the subsequent formulas simpler, we introduce the notation
\begin{equation}\label{eq40}
\ell \left( \theta  \right) = \begin{cases} \left|\csc\left(2\theta\right)\right| & \text{ if } \; \frac{\pi}{4} < \left|\theta\right| <\frac{\pi}{2} \\ 1 & \text{ if } \; \left|\theta\right| \leq \frac{\pi}{4} . \end{cases}
\end{equation}
If we substitute the integral representation (see \cite[10.32.E9]{NIST})
\[
K_\nu  \left( t \right) = \int_0^{ + \infty } {e^{ - t\cosh s} \cosh \left( {\nu s} \right)ds} 
\]
into \eqref{eq2} and perform the change of variable $u = t\cosh s$, we find that
\begin{equation}\label{eq30}
R_N \left( {z,\mu ,\nu } \right) = \left( { - 1} \right)^N \frac{{2^{\mu  + 1} z^{\mu  - 2N - 1} }}{{\Gamma \left( {\frac{{\nu  - \mu  + 1}}{2}} \right)\Gamma \left( {\frac{{1 - \mu  - \nu }}{2}} \right)}}\int_0^{ + \infty } {u^{2N - \mu } e^{ - u} \int_0^{ + \infty } {\frac{{\cosh \left( {\nu s} \right)\cosh ^{\mu - 2N - 1} s}}{{1 + \left( {u/z\cosh s} \right)^2 }}ds} du} 
\end{equation}
for any non-negative integer $N$ with $\left|\arg z\right| < \frac{\pi}{2}$ and $\left| {\Re \left( \nu  \right)} \right| + \Re \left( \mu  \right) < 2N + 1$. Using the equality
\[
a_N \left( { - \mu ,\nu } \right) = \left( { - 1} \right)^N z^{ - \mu  + 2N + 1} \left( {R_N \left( {z,\mu ,\nu } \right) - R_{N + 1} \left( {z,\mu ,\nu } \right)} \right),
\]
we deduce that
\begin{equation}\label{eq33}
a_N \left( { - \mu ,\nu } \right) = \frac{{2^{\mu  + 1} }}{{\Gamma \left( {\frac{{\nu  - \mu  + 1}}{2}} \right)\Gamma \left( {\frac{{1 - \mu  - \nu }}{2}} \right)}}\int_0^{ + \infty } {u^{2N - \mu } e^{ - u} \int_0^{ + \infty } {\cosh \left( {\nu s} \right)\cosh ^{\mu  - 2N - 1} sds} du} 
\end{equation}
as long as $\left| {\Re \left( \nu  \right)} \right| + \Re \left( \mu  \right) < 2N + 1$. The transformations we have made are not necessary for the derivation of the bounds for the right-half plane, but will be important when we extend them beyond the imaginary axis. Noting that for any $r>0$, $1/\left| {1 + re^{ - 2\theta i} } \right| \le \ell \left( \theta  \right)$, trivial estimation of \eqref{eq30} and the formula \eqref{eq33} provide the error bound
\begin{equation}\label{eq31}
\left| {R_N \left( {z,\mu ,\nu } \right)} \right| \le  \left|\frac{{\Gamma \left( {\frac{{\Re \left( \nu  \right) - \Re \left( \mu  \right) + 1}}{2}} \right)\Gamma \left( {\frac{{1 - \Re \left( \mu  \right) - \Re \left( \nu  \right)}}{2}} \right)}}{\Gamma \left( {\frac{{\nu  - \mu  + 1}}{2}} \right)\Gamma \left( {\frac{{1 - \mu  - \nu }}{2}} \right)}\right| \left| {z^{\mu  - 1} } \right|\frac{\left|a_N \left( { - \Re \left( \mu  \right),\Re \left( \nu  \right)} \right)\right|}{{\left| z \right|^{2N} }} \ell \left( \theta  \right),
\end{equation}
using the notation $\theta = \arg z$. If $\Re \left( \nu  \right) + \Re \left( \mu  \right) \neq \nu + \mu$ and $\Re \left( \nu  \right) - \Re \left( \mu  \right) \neq \nu - \mu$, and neither number on the right-hand side is a negative odd integer, then we can use the inequality
\begin{gather}\label{eq36}
\begin{split}
\left|\frac{{\Gamma \left( {\frac{{\Re \left( \nu  \right) - \Re \left( \mu  \right) + 1}}{2}} \right)\Gamma \left( {\frac{{1 - \Re \left( \mu  \right) - \Re \left( \nu  \right)}}{2}} \right)}}{\Gamma \left( {\frac{{\nu  - \mu  + 1}}{2}} \right)\Gamma \left( {\frac{{1 - \mu  - \nu }}{2}} \right)}\right|  & \le \left| {\frac{{\cos \left( {\pi \frac{{\nu  - \mu }}{2}} \right)\cos \left( {\pi \frac{{\mu  + \nu }}{2}} \right)}}{{\cos \left( {\pi \frac{{\Re \left( \nu  \right) - \Re \left( \mu  \right)}}{2}} \right)\cos \left( {\pi \frac{{\Re \left( \mu  \right) + \Re \left( \nu  \right)}}{2}} \right)}}} \right| \\
& = \left| {\frac{{\cos \left( {\pi \mu } \right) + \cos \left( {\pi \nu } \right)}}{{\cos \left( {\pi \Re \left( \mu  \right)} \right) + \cos \left( {\pi \Re \left( \nu  \right)} \right)}}} \right|,
\end{split}
\end{gather}
which follows from the reflection formula for the Gamma function and the inequality $\left| {\Gamma \left( w \right)} \right| \le \left| {\Gamma \left( {\Re \left( w \right)} \right)} \right|$, to simplify the error bound \eqref{eq31} to
\begin{equation}\label{eq32}
\left| {R_N \left( {z,\mu ,\nu } \right)} \right| \le \left| {\frac{{\cos \left( {\pi \mu } \right) + \cos \left( {\pi \nu } \right)}}{{\cos \left( {\pi \Re \left( \mu  \right)} \right) + \cos \left( {\pi \Re \left( \nu  \right)} \right)}}} \right|\left| {z^{\mu  - 1} } \right|\frac{\left|a_N \left( { - \Re \left( \mu  \right),\Re \left( \nu  \right)} \right)\right|}{{\left| z \right|^{2N} }} \ell \left( \theta  \right).
\end{equation}
If $\Re \left( \nu  \right) + \Re \left( \mu  \right) \neq \nu + \mu$ or $\Re \left( \nu  \right) - \Re \left( \mu  \right) \neq \nu - \mu$ and at least one of them is an positive odd integer, then the limiting value has to be taken in \eqref{eq32}. The existence of the limit follows from the definition of the coefficients $a_N \left( { - \Re \left( \mu  \right),\Re \left( \nu  \right)} \right)$ and the fact that $\left| {\Re \left( \nu  \right)} \right| + \Re \left( \mu  \right) < 2N + 1$.

From \eqref{eq32} (or \eqref{eq31}) it is seen that when both $\mu$ and $\nu$ are real, the absolute value of the remainder term $R_N \left( {z,\mu ,\nu } \right)$ is bounded by the absolute value of the first omitted term of the asymptotic series multiplied by $\ell \left( \theta  \right)$, provided that $\left| \nu  \right| + \mu  < 2N + 1$. In addition, if $z>0$ we have $0 < 1/\left( {1 + \left( {u/z\cosh s} \right)^2 } \right) < 1$ in \eqref{eq30}, and the mean value theorem of integration shows that
\[
R_N \left( {z,\mu ,\nu } \right) = z^{\mu  - 1} \left( { - 1} \right)^N \frac{{a_N \left( { - \mu ,\nu } \right)}}{z^{2N}}\Theta ,
\]
as long as $\left| \nu  \right| + \mu  < 2N + 1$, where $0<\Theta<1$ is an appropriate number depending on $z$, $\mu$, $\nu$ and $N$.

Our bounds for $R_N \left( {z,\mu ,\nu } \right)$ are unrealistic near the Stokes lines $\arg z = \theta = \pm \frac{\pi }{2}$ due to the presence of the factor $\left|\csc \left(2\theta\right)\right|$. A better bound for $R_N \left( {z,\mu ,\nu } \right)$ near these lines can be derived as follows. Let $0 < \varphi  < \frac{\pi }{2}$ be an acute angle that may depend on $\mu$ and $N$ and suppose that $\frac{\pi }{4} + \varphi  < \theta  < \frac{\pi }{2} + \varphi$. An analytic continuation of the representation \eqref{eq30} to this sector can be found by rotating the path of integration through the angle $\varphi$, to obtain
\begin{align*}
& R_N \left( {z,\mu ,\nu } \right) = \left( { - 1} \right)^N \frac{{2^{\mu  + 1} z^{\mu  - 2N - 1} }}{{\Gamma \left( {\frac{{\nu  - \mu  + 1}}{2}} \right)\Gamma \left( {\frac{{1 - \mu  - \nu }}{2}} \right)}}\int_0^{ + \infty e^{i\varphi } } {u^{2N - \mu } e^{ - u} \int_0^{ + \infty } {\frac{{\cosh \left( {\nu s} \right)\cosh ^{\mu  - 2N - 1} s}}{{1 + \left( {u/z\cosh s} \right)^2 }}ds} du} 
\\ & = \left( { - 1} \right)^N \frac{{2^{\mu  + 1} z^{\mu  - 2N - 1} }}{{\Gamma \left( {\frac{{\nu  - \mu  + 1}}{2}} \right)\Gamma \left( {\frac{{1 - \mu  - \nu }}{2}} \right)}}\left( {\frac{{e^{i\varphi } }}{{\cos \varphi }}} \right)^{2N + 1 - \mu } \int_0^{ + \infty } {t^{2N - \mu } e^{ - \frac{{te^{i\varphi } }}{{\cos \varphi }}} \int_0^{ + \infty } {\frac{{\cosh \left( {\nu s} \right)\cosh ^{\mu  - 2N - 1} s}}{{1 + \left( {te^{i\varphi } /z\cosh s\cos \varphi } \right)^2 }}ds} dt} .
\end{align*}
Employing the inequality $1/\left| {1 + re^{ - 2\theta i} } \right| \le \ell \left( \theta  \right)$ ($r>0$) and the expression \eqref{eq33}, we deduce the error bound
\begin{gather}\label{eq35}
\begin{split}
& \left| {R_N \left( {z,\mu ,\nu } \right)} \right| \le \\ &
\frac{{e^{\Im \left( \mu  \right)\varphi } \csc \left( {2\left( {\theta  - \varphi } \right)} \right)}}{{\left( {\cos \varphi } \right)^{2N + 1 - \Re \left( \mu  \right)} }}\left| {\frac{{\Gamma \left( {\frac{{\Re \left( \nu  \right) - \Re \left( \mu  \right) + 1}}{2}} \right)\Gamma \left( {\frac{{1 - \Re \left( \mu  \right) - \Re \left( \nu  \right)}}{2}} \right)}}{{\Gamma \left( {\frac{{\nu  - \mu  + 1}}{2}} \right)\Gamma \left( {\frac{{1 - \mu  - \nu }}{2}} \right)}}} \right|\left| {z^{\mu  - 1} } \right|\frac{{\left| {a_N \left( { - \Re \left( \mu  \right),\Re \left( \nu  \right)} \right)} \right|}}{{\left| z \right|^{2N} }},
\end{split}
\end{gather}
provided that $\frac{\pi }{4} < \frac{\pi }{4} + \varphi  < \theta  < \frac{\pi }{2} + \varphi < \pi$ and $\left| {\Re \left( \nu  \right)} \right| + \Re \left( \mu  \right) < 2N + 1$. It seems that we can not minimise the factor in the front as a function of $\varphi$ in simple terms. Nevertheless, if we concentrate on the function
\begin{equation}\label{eq34}
\frac{ \csc \left( {2\left( {\theta  - \varphi } \right)} \right)}{\left( {\cos \varphi } \right)^{2N + 1 - \Re \left( \mu  \right)}},
\end{equation}
the minimisation can be done by applying a lemma of Meijer's \cite[p. 956]{Meijer}. In our case, Meijer's lemma gives that the minimising value $\varphi = \varphi^\ast$ in \eqref{eq34}, is the unique solution of the implicit equation
\[
\left( {2N + 3 - \Re \left( \mu  \right)} \right)\cos \left( {3\varphi^\ast - 2\theta} \right) = \left( {2N - 1 - \Re \left( \mu  \right)} \right)\cos \left( {\varphi^\ast- 2\theta } \right),
\]
that satisfies $-\frac{\pi}{2}+\theta < \varphi^\ast <\frac{\pi}{2}$ if $\frac{3\pi}{4} \leq \theta <\pi$; $-\frac{\pi}{2}+\theta < \varphi^\ast <-\frac{\pi}{4}+\theta$ if $\frac{\pi}{2} \leq \theta <\frac{3\pi}{4}$; and $0 < \varphi^\ast <-\frac{\pi}{4}+\theta$ if $\frac{\pi}{4} < \theta <\frac{\pi}{2}$. With this choice of $\varphi$, \eqref{eq35} provides an error bound for the range $\frac{\pi}{4}<\arg z <\pi$. Again, if $\Re \left( \nu  \right) + \Re \left( \mu  \right) \neq \nu + \mu$ and $\Re \left( \nu  \right) - \Re \left( \mu  \right) \neq \nu - \mu$, and neither number on the right-hand side is a negative odd integer, then we can employ the inequality \eqref{eq36} to simplify \eqref{eq35}.

We can make our bounds simpler if $\arg z$ is close to $\frac{\pi}{2}$ as follows. When $\arg z = \theta = \frac{\pi}{2}$, the minimising value $\varphi^\ast$ in \eqref{eq34} is given explicitly by
\[
\varphi^\ast  = \arctan \left( {\frac{1}{{\sqrt {2N + 2 - \Re \left( \mu  \right)} }}} \right),
\]
and therefore we have
\begin{align*}
\frac{{\csc \left( {2\left( {\theta  - \varphi^\ast } \right)} \right)}}{{\left( {\cos \varphi^\ast } \right)^{2N + 1 - \Re \left( \mu  \right)} }} \le \frac{{\csc \left( {2\left( {\frac{\pi }{2} - \varphi^\ast } \right)} \right)}}{{\left( {\cos \varphi^\ast } \right)^{2N + 1 - \Re \left( \mu  \right)} }} & = \frac{1}{2}\left( {1 + \frac{1}{{2N + 2 - \Re \left( \mu  \right)}}} \right)^{N + \frac{{3-\Re \left( \mu  \right)}}{2}} \sqrt {2N + 2 - \Re \left( \mu  \right)} \\ & \le \frac{1}{2}\sqrt {e\left( {2N + \frac{5}{2} - \Re \left( \mu  \right)} \right)}, 
\end{align*}
as long as $\frac{\pi }{4} + \varphi^\ast < \theta  \leq \frac{\pi }{2}$, $N\geq 0$ and $\left| {\Re \left( \nu  \right)} \right| + \Re \left( \mu  \right)< 2N + 1$.

The corresponding error bounds for the range $-\pi<\arg z <-\frac{\pi}{4}$ can be obtained from these results and the fact that $\left| {R_N \left( {z,\mu ,\nu } \right)} \right| = \left| {\overline {R_N \left( {z,\mu ,\nu } \right)} } \right| = \left| {R_N \left( {\bar z,\bar \mu ,\bar \nu } \right)} \right|$.

\section{Exponentially improved asymptotic expansions}\label{section4}

We shall find it convenient to express our exponentially improved expansion in terms of the (scaled) Terminant function, which is defined in terms of the Incomplete gamma function as
\[
\widehat T_p \left( w \right) = \frac{{e^{\pi ip} \Gamma \left( p \right)}}{{2\pi i}}\Gamma \left( {1 - p,w} \right) = \frac{e^{\pi ip} w^{1 - p} e^{ - w} }{2\pi i}\int_0^{ + \infty } {\frac{{t^{p - 1} e^{ - t} }}{w + t}dt} \; \text{ for } \; \Re\left(p\right)>0 \; \text{ and } \; \left| \arg w \right| < \pi ,
\]
and by analytic continuation elsewhere. Olver \cite[equations (4.5) and (4.6)]{Olver4} showed that when $\left|p\right| \sim \left|w\right|$ and $w \to \infty$, we have
\begin{equation}\label{eq11}
ie^{ - \pi ip} \widehat T_p \left( w \right) = \begin{cases} \mathcal{O}\left( {e^{ - w - \left| w \right|} } \right) & \; \text{ if } \; \left| {\arg w} \right| \le \pi \\ \mathcal{O}\left(1\right) & \; \text{ if } \; - 3\pi  < \arg w \le  - \pi. \end{cases}
\end{equation}
Concerning the smooth transition of the Stokes discontinuities, we will use the more precise asymptotic formulas
\begin{equation}\label{eq12}
\widehat T_p \left( w \right) = \frac{1}{2} + \frac{1}{2}\mathop{\text{erf}} \left( {c\left( \varphi  \right)\sqrt {\frac{1}{2}\left| w \right|} } \right) + \mathcal{O}\left( {\frac{{e^{ - \frac{1}{2}\left| w \right|c^2 \left( \varphi  \right)} }}{{\left| w \right|^{\frac{1}{2}} }}} \right)
\end{equation}
for $-\pi +\delta \leq \arg w \leq 3 \pi -\delta$, $0 < \delta  \le 2\pi$; and
\begin{equation}\label{eq13}
e^{ - 2\pi ip} \widehat T_p \left( w \right) =  - \frac{1}{2} + \frac{1}{2}\mathop{\text{erf}} \left( { - \overline {c\left( { - \varphi } \right)} \sqrt {\frac{1}{2}\left| w \right|} } \right) + \mathcal{O}\left( {\frac{{e^{ - \frac{1}{2}\left| w \right|\overline {c^2 \left( { - \varphi } \right)} } }}{{\left| w \right|^{\frac{1}{2}} }}} \right)
\end{equation}
for $- 3\pi  + \delta  \le \arg w \le \pi  - \delta$, $0 < \delta \le 2\pi$. Here $\varphi = \arg w$ and $\mathop{\text{erf}}$ denotes the Error function. The quantity $c\left( \varphi  \right)$ is defined implicitly by the equation
\[
\frac{1}{2}c^2 \left( \varphi  \right) = 1 + i\left( {\varphi  - \pi } \right) - e^{i\left( {\varphi  - \pi } \right)},
\]
and corresponds to the branch of $c\left( \varphi  \right)$ which has the following expansion in the neighbourhood of $\varphi = \pi$:
\begin{equation}\label{eq14}
c\left( \varphi  \right) = \left( {\varphi  - \pi } \right) + \frac{i}{6}\left( {\varphi  - \pi } \right)^2  - \frac{1}{{36}}\left( {\varphi  - \pi } \right)^3  - \frac{i}{{270}}\left( {\varphi  - \pi } \right)^4  +  \cdots .
\end{equation}
For complete asymptotic expansions, see Olver \cite{Olver5}. We remark that Olver uses the different notation $F_p \left( w \right) = ie^{ - \pi ip} \widehat T_p \left( w \right)$ for the Terminant function and the other branch of the function $c\left( \varphi  \right)$. For further properties of the Terminant function, see, for example, Paris and Kaminski \cite[Chapter 6]{Paris3}.

\subsection{Proof of Theorems \ref{thm2} and \ref{thm3}} First, we derive order estimates for the remainder $R_{N,M} \left( {z,\mu ,\nu } \right)$ with subject to the condition that $\left| {\Re \left( \nu  \right)} \right| < M + \frac{1}{2}$. First, we suppose that $\left|\arg z\right|< \frac{\pi}{2}$. Let $M$ be a fixed non-negative integer. Substituting the expression \eqref{eq28} into \eqref{eq2} and using the definition of the Terminant function we find that
\begin{multline*}
R_N \left( {z,\mu ,\nu } \right) = \frac{{2^{\mu  + 1} }}{{\Gamma \left( {\frac{{\nu  - \mu  + 1}}{2}} \right)\Gamma \left( {\frac{{1 - \mu  - \nu }}{2}} \right)}}\left( {\pi e^{\frac{\pi }{2}i\mu } \left( {\frac{\pi}{-2iz}} \right)^{\frac{1}{2}} e^{iz} \sum\limits_{m = 0}^{M - 1} {\frac{{a_m \left( \nu  \right)}}{{\left( {-iz} \right)^m }} \widehat T_{2N - m - \mu  + \frac{1}{2}} \left( {iz} \right)} }\right. \\ \left.{+ \pi e^{ - \frac{\pi }{2}i\mu } \left( {\frac{\pi }{{2iz }}} \right)^{\frac{1}{2}} e^{ - iz} \sum\limits_{m = 0}^{M - 1} {\frac{{a_m \left( \nu  \right)}}{{\left( {iz } \right)^m }}e^{2\pi i\mu } \widehat T_{2N - m - \mu  + \frac{1}{2}} \left( { - iz} \right)}  + R_{N,M} \left( {z,\mu ,\nu } \right)} \right),
\end{multline*}
with
\begin{gather}\label{eq9}
\begin{split}
R_{N,M} \left( {z,\mu ,\nu } \right) & = \left( { - 1} \right)^N \frac{{z^{\mu  - 2N - 1} }}{2}\left( {\int_0^{ + \infty } {\frac{{t^{2N - \mu } K_M \left( {t,\nu } \right)}}{{1 - it/z}}dt}  + \int_0^{ + \infty } {\frac{{t^{2N - \mu } K_M \left( {t,\nu } \right)}}{{1 + it/z}}dt} } \right)
\\ & = \left( { - 1} \right)^N \frac{{\left( {e^{i\theta } } \right)^{\mu  - 2N - 1} }}{2}\left( {\int_0^{ + \infty } {\frac{{\tau ^{2N - \mu } K_M \left( {r\tau ,\nu } \right)}}{{1 - i\tau e^{ - i\theta } }}d\tau }  + \int_0^{ + \infty } {\frac{{\tau ^{2N - \mu } K_M \left( {r\tau ,\nu } \right)}}{{1 + i\tau e^{ - i\theta } }}d\tau } } \right),
\end{split}
\end{gather}
under the assumption that $2N - M + \frac{1}{2} > \Re \left( \mu  \right)$. Here we have taken $z = re^{i\theta }$. Using the integral formula \eqref{eq29}, $K_M \left( {r\tau ,\nu } \right)$ can be written as
\begin{multline*}
K_M \left( {r\tau ,\nu } \right) = \left( { - 1} \right)^M \frac{{\cos \left( {\pi \nu } \right)}}{\pi }\left( {r\tau } \right)^{ - M - \frac{1}{2}} e^{ - r\tau } \int_0^{ + \infty } {\frac{{s^{M - \frac{1}{2}} e^{ - s} K_\nu  \left( s \right)}}{{1 + s/r}}ds} \\ + \left( { - 1} \right)^M \frac{{\cos \left( {\pi \nu } \right)}}{\pi }\left( {r\tau } \right)^{ - M - \frac{1}{2}} e^{ - r\tau } \left( {\tau  - 1} \right)\int_0^{ + \infty } {\frac{{s^{M - \frac{1}{2}} e^{ - s} K_\nu  \left( s \right)}}{{\left( {1 + r\tau /s} \right)\left( {1 + s/r} \right)}}ds} ,
\end{multline*}
provided that $\left| {\Re \left( \nu  \right)} \right| < M + \frac{1}{2}$. Noting that
\[
0 < \frac{1}{{1 + s/r}},\frac{1}{{\left( 1 + r\tau /s \right)\left( {1 + s/r} \right)}} < 1
\]
for positive $r$, $\tau$ and $s$, substitution into \eqref{eq9} yields the upper bound
\begin{align*}
& \left| {R_{N,M} \left( {z,\mu ,\nu } \right)} \right| \le  \pi \left| {e^{\frac{\pi }{2}i\mu } } \right|\left( {\frac{\pi }{{2\left| z \right|}}} \right)^{\frac{1}{2}} \frac{{\left| {\cos \left( {\pi \nu } \right)} \right|}}{{\left| {\cos \left( {\pi \Re \left( \nu  \right)} \right)} \right|}}\frac{{\left| {a_M \left( \Re\left(\nu \right) \right)} \right|}}{{\left| z \right|^M }}\left| {e^{i\theta \mu } \frac{{e^{ - \frac{\pi }{2}i\mu } }}{{2\pi }}\int_0^{ + \infty } {\frac{{\tau ^{2N - M - \mu  - \frac{1}{2}} e^{ - r\tau } }}{{1 - i\tau e^{ - i\theta } }}d\tau } } \right|\\
& +  \left( {\frac{\pi }{{2\left| z \right|}}} \right)^{\frac{1}{2}} \frac{{\left| {\cos \left( {\pi \nu } \right)} \right|}}{{\left| {\cos \left( {\pi \Re \left( \nu  \right)} \right)} \right|}}\frac{{\left| {a_M \left( \Re\left(\nu \right) \right)} \right|}}{{\left| z \right|^M }}\left| {e^{i\theta \mu } } \right|\frac{1}{{2 }}\int_0^{ + \infty } {\tau ^{2N - M - \Re \left( \mu  \right) - \frac{1}{2}} e^{ - r\tau } \left| {\frac{{\tau  - 1}}{{\tau  + ie^{i\theta } }}} \right|d\tau } \\
& + \pi \left| {e^{ - \frac{\pi }{2}i\mu } } \right|\left( {\frac{\pi }{{2\left| z \right|}}} \right)^{\frac{1}{2}} \frac{{\left| {\cos \left( {\pi \nu } \right)} \right|}}{{\left| {\cos \left( {\pi \Re \left( \nu  \right)} \right)} \right|}}\frac{{\left| {a_M \left( \Re\left(\nu \right)  \right)} \right|}}{{\left| z \right|^M }}\left| {e^{i\theta \mu } \frac{e^{\frac{\pi }{2}i\mu }}{{2\pi }}\int_0^{ + \infty } {\frac{{\tau ^{2N - M - \mu  - \frac{1}{2}} e^{ - r\tau } }}{{1 + i\tau e^{ - i\theta } }}d\tau } } \right| \\
& + \left( {\frac{\pi }{{2\left| z \right|}}} \right)^{\frac{1}{2}} \frac{{\left| {\cos \left( {\pi \nu } \right)} \right|}}{{\left| {\cos \left( {\pi \Re \left( \nu  \right)} \right)} \right|}}\frac{{\left| {a_M \left( \Re\left(\nu \right)  \right)} \right|}}{{\left| z \right|^M }}\left| {e^{i\theta \mu }    } \right|\frac{1}{{2 }}\int_0^{ + \infty } {\tau ^{2N - M - \Re \left( \mu  \right) - \frac{1}{2}} e^{ - r\tau } \left| {\frac{{\tau  - 1}}{{\tau  - ie^{i\theta } }}} \right|d\tau } .
\end{align*}
Since $\left| {\left( {\tau  - 1} \right)/\left( {\tau  \pm ie^{i\theta } } \right)} \right| \le 1$, we find that
\begin{align*}
\left| {R_{N,M} \left( {z,\mu ,\nu } \right)} \right| \le \; & \pi \left| {e^{\frac{\pi }{2}i\mu } } \right|\left( {\frac{\pi }{{2\left| z \right|}}} \right)^{\frac{1}{2}} \frac{{\left| {\cos \left( {\pi \nu } \right)} \right|}}{{\left| {\cos \left( {\pi \Re \left( \nu  \right)} \right)} \right|}}\frac{{\left| {a_M \left( \Re\left(\nu \right)  \right)} \right|}}{{\left| z \right|^M }}\left| {e^{iz} \widehat T_{2N - M - \mu  + \frac{1}{2}} \left( {iz} \right)} \right| \\ & + \pi \left| {e^{ - \frac{\pi }{2}i\mu } } \right|\left( {\frac{\pi }{{2\left| z \right|}}} \right)^{\frac{1}{2}} \frac{{\left| {\cos \left( {\pi \nu } \right)} \right|}}{{\left| {\cos \left( {\pi \Re \left( \nu  \right)} \right)} \right|}}\frac{{\left| {a_M \left( \Re\left(\nu \right) \right)} \right|}}{{\left| z \right|^M }}\left| {e^{2\pi i\mu } e^{ - iz} \widehat T_{2N - M - \mu  + \frac{1}{2}} \left( { - iz} \right)} \right|
\\ & + \left( {\frac{\pi }{2}} \right)^{\frac{1}{2}} \left| {z^\mu  } \right|\frac{{\left| {\cos \left( {\pi \nu } \right)} \right|}}{{\left| {\cos \left( {\pi \Re \left( \nu  \right)} \right)} \right|}}\frac{{\left| {a_M \left( \Re\left(\nu \right) \right)} \right|\Gamma \left( {2N - M - \Re \left( \mu  \right) + \frac{1}{2}} \right)}}{{\left| z \right|^{2N + 1} }}.
\end{align*}
By continuity, this bound holds in the closed sector $\left|\arg z\right| \leq \frac{\pi}{2}$. Assume that $N = \frac{1}{2}\left| z \right| + \rho$ where $\rho$ is bounded. Employing Stirling's formula, we find that
\[
\left( {\frac{\pi }{2}} \right)^{\frac{1}{2}} \left| {z^\mu  } \right|\frac{{\left| {\cos \left( {\pi \nu } \right)} \right|}}{{\left| {\cos \left( {\pi \Re \left( \nu  \right)} \right)} \right|}}\frac{{\left| {a_M \left( \Re \left( \nu  \right) \right)} \right|\Gamma \left( {2N - M - \Re \left( \mu  \right) + \frac{1}{2}} \right)}}{{\left| z \right|^{2N + 1} }} = \mathcal{O}_{M,\mu ,\rho } \left( {\frac{{e^{ - \left| z \right|} }}{{\left| z \right|}}\frac{{\left| {\cos \left( {\pi \nu } \right)} \right|}}{{\left| {\cos \left( {\pi \Re \left( \nu  \right)} \right)} \right|}}\frac{{\left| {a_M \left(\Re \left( \nu  \right)\right)} \right|}}{{\left| z \right|^M }}} \right)
\]
as $z \to \infty$. Olver's estimation \eqref{eq11} shows that
\[
\left| {e^{ \pm iz} \widehat T_{2N - M - \mu  + \frac{1}{2}} \left( { \pm iz} \right)} \right| = \mathcal{O}_{M,\mu } \left( { e^{ - \left| z \right|} } \right)
\]
for large $z$. Therefore, we have that
\begin{equation}\label{eq10}
R_{N,M} \left( {z,\mu ,\nu } \right) = \mathcal{O}_{M,\mu ,\rho } \left( {\frac{e^{ - \left| z \right|} }{\left| z \right|^{\frac{1}{2}}}\frac{{\left| {\cos \left( {\pi \nu } \right)} \right|}}{{\left| {\cos \left( {\pi \Re \left( \nu  \right)} \right)} \right|}}\frac{{\left| {a_M \left(\Re \left( \nu  \right)\right)} \right|}}{{\left| z \right|^M }}} \right)
\end{equation}
as $z \to \infty$ in the sector $\left|\arg z\right| \leq \frac{\pi}{2}$.

Consider now the sectors $\frac{\pi }{2} < \pm \arg z < \frac{{3\pi }}{2}$. When $z$ enters the sector $\frac{\pi}{2}< \arg z < \frac{3\pi}{2}$, the pole in the first integral in \eqref{eq9} crosses the integration path. Similarly, when $z$ enters the sector $-\frac{3\pi}{2}< \arg z < -\frac{\pi}{2}$, the pole in the second integral in \eqref{eq9} crosses the integration path. According to the residue theorem, we obtain
\begin{gather}\label{eq16}
\begin{split}
R_{N,M} \left( {z,\mu ,\nu } \right) = \; & \pi e^{ \pm \frac{\pi }{2}i\mu } K_M \left( {ze^{ \mp \frac{\pi }{2}i} ,\nu } \right) \\ & + \left( { - 1} \right)^N \frac{{z^{\mu  - 2N - 1} }}{2}\left( {\int_0^{ + \infty } {\frac{{t^{2N - \mu } K_M \left( {t,\nu } \right)}}{{1 - it/z}}dt}  + \int_0^{ + \infty } {\frac{{t^{2N - \mu } K_M \left( {t,\nu } \right)}}{{1 + it/z}}dt} } \right)\\
 = \; & \pi e^{ \pm \frac{\pi }{2}i\mu } K_M \left( {ze^{ \mp \frac{\pi }{2}i} ,\nu } \right) - e^{ \pm \pi i\mu } R_{N,M} \left( {ze^{ \mp \pi i} ,\mu ,\nu } \right)
\end{split}
\end{gather}
when $\frac{\pi }{2} <  \pm \arg z < \frac{{3\pi }}{2}$. It follows that
\[
\left| {R_{N,M} \left( {z,\mu ,\nu } \right)} \right| \le \pi e^{ \mp \frac{\pi }{2}\Im \left( \mu  \right)} \left| {K_M \left( {ze^{ \mp \frac{\pi }{2}i} ,\nu } \right)} \right| + e^{ \mp \pi \Im \left( \mu  \right)} \left| {R_{N,M} \left( {ze^{ \mp \pi i} ,\mu ,\nu } \right)} \right|
\]
in the closed sectors $\frac{\pi }{2} \leq \pm \arg z \leq \frac{{3\pi }}{2}$, using continuity. Meijer \cite{Meijer} proved that
\[
K_M \left( {ze^{ \mp \frac{\pi }{2}i} ,\nu } \right) = \mathcal{O}_M \left( {\frac{{e^{ \mp \Im \left( z \right)} }}{{\left| z \right|^{\frac{1}{2}} }}\frac{{\left| {\cos \left( {\pi \nu } \right)} \right|}}{{\left| {\cos \left( {\pi \Re \left( \nu  \right)} \right)} \right|}}\frac{{\left| {a_M \left( \Re \left( \nu  \right)  \right)} \right|}}{{\left| z \right|^M }}} \right)
\]
if $\left| {\arg \left( {ze^{ \mp \frac{\pi }{2}i} } \right)} \right| \leq \pi$ and $\left| {\Re \left( \nu  \right)} \right| < M + \frac{1}{2}$. Combining this bound together with \eqref{eq10}, yields
\[
R_{N,M} \left( {z,\mu ,\nu } \right) = \mathcal{O}_{M,\mu ,\rho } \left( {\frac{{e^{ \mp \Im \left( z \right)} }}{{\left| z \right|^{\frac{1}{2}} }}\frac{{\left| {\cos \left( {\pi \nu } \right)} \right|}}{{\left| {\cos \left( {\pi \Re \left( \nu  \right)} \right)} \right|}}\frac{{\left| {a_M \left( \Re \left( \nu  \right) \right)} \right|}}{{\left| z \right|^M }}} \right)
\]
as $z \to \infty$ in the sector $\frac{\pi }{2} \leq \pm \arg z \leq \frac{{3\pi }}{2}$.

Now, let $M$ be an arbitrary fixed non-negative integer, and let $M'$ be a non-negative integer such that $\left| {\Re \left( \nu  \right)} \right| < M' + \frac{1}{2}$. We have
\begin{multline*}
R_{N,M} \left( {z,\mu ,\nu } \right) = \pi e^{\frac{\pi }{2}i\mu } \left( {\frac{\pi }{{ - 2iz}}} \right)^{\frac{1}{2}} e^{iz} \sum\limits_{m = M}^{M' - 1} {\frac{{a_m \left( \nu  \right)}}{{\left( { - iz} \right)^m }}\widehat T_{2N - m - \mu  + \frac{1}{2}} \left( {iz} \right)} 
\\ + \pi e^{ - \frac{\pi }{2}i\mu } \left( {\frac{\pi }{{2iz}}} \right)^{\frac{1}{2}} e^{ - iz} \sum\limits_{m = M}^{M' - 1} {\frac{{a_m \left( \nu  \right)}}{{\left( {iz} \right)^m }}e^{2\pi i\mu } \widehat T_{2N - m - \mu  + \frac{1}{2}} \left( { - iz} \right)}  + R_{N,M'} \left( {z,\mu ,\nu } \right).
\end{multline*}
Trivial estimation yields
\begin{multline*}
\left| {R_{N,M} \left( {z,\mu ,\nu } \right)} \right| \le \pi \left| {e^{\frac{\pi }{2}i\mu } } \right|\left( {\frac{\pi }{{2\left| z \right|}}} \right)^{\frac{1}{2}} \sum\limits_{m = M}^{M' - 1} {\frac{{\left| {a_m \left( \nu  \right)} \right|}}{{\left| z \right|^m }}\left| {e^{iz} \widehat T_{2N - m - \mu  + \frac{1}{2}} \left( {iz} \right)} \right|} 
\\  + \pi \left| {e^{ - \frac{\pi }{2}i\mu } } \right|\left( {\frac{\pi }{{2\left| z \right|}}} \right)^{\frac{1}{2}} \sum\limits_{m = M}^{M' - 1} {\frac{{\left| {a_m \left( \nu  \right)} \right|}}{{\left| z \right|^m }}\left| {e^{2\pi i\mu } } \right|\left| {e^{ - iz} \widehat T_{2N - m - \mu  + \frac{1}{2}} \left( { - iz} \right)} \right|}  + \left| {R_{N,M'} \left( {z,\mu ,\nu } \right)} \right|.
\end{multline*}
Employing the previously obtained bounds for $R_{N,M'} \left( {z,\mu ,\nu } \right)$ and Olver's estimation \eqref{eq11} together with the connection formula for the Terminant function \cite[p. 260]{Paris3}, shows that $R_{N,M} \left( {z,\mu ,\nu } \right)$ indeed satisfies the order estimates prescribed in Theorem \ref{thm2}.

\subsection{Stokes phenomenon and Berry's transition} We study the Stokes phenomenon related to the large-$z$ asymptotic expansion of $S_{\mu,\nu}\left(z\right)$ occurring when $\arg z$ passes through the values $\pm \frac{\pi}{2}$. In the range $\left|\arg z\right|<\frac{\pi}{2}$, the asymptotic expansion
\begin{equation}\label{eq15}
S_{\mu ,\nu } \left( z \right) \sim z^{\mu  - 1} \sum\limits_{n = 0}^\infty  {\left( { - 1} \right)^n \frac{{a_n \left( { - \mu ,\nu } \right)}}{{z^{2n} }}} 
\end{equation}
holds as $z \to \infty$. From \eqref{eq16} we have
\begin{align*}
S_{\mu ,\nu } \left( z \right) & = \frac{{2^{\mu  + 1} }}{{\Gamma \left( {\frac{{\nu  - \mu  + 1}}{2}} \right)\Gamma \left( {\frac{{1 - \mu  - \nu }}{2}} \right)}} R_{0,0} \left( {z,\mu ,\nu } \right) \\ & = \frac{{2^{\mu  + 1} }}{{\Gamma \left( {\frac{{\nu  - \mu  + 1}}{2}} \right)\Gamma \left( {\frac{{1 - \mu  - \nu }}{2}} \right)}} \left(\pi e^{\frac{\pi }{2}i\mu } K_0 \left( {ze^{ - \frac{\pi }{2}i} ,\nu } \right) -  e^{\pi i\mu } R_{0,0} \left( {ze^{ - \pi i} ,\mu ,\nu } \right) \right)\\ & = \frac{{2^{\mu  + 1} }}{{\Gamma \left( {\frac{{\nu  - \mu  + 1}}{2}} \right)\Gamma \left( {\frac{{1 - \mu  - \nu }}{2}} \right)}} \pi e^{\frac{\pi }{2}i\mu } K_\nu  \left( { - iz} \right) - e^{\pi i\mu } S_{\mu ,\nu } \left( {ze^{ - \pi i} } \right)
\end{align*}
when $\frac{\pi}{2} < \arg z < \frac{3\pi}{2}$. Similarly, from \eqref{eq16} we find
\[
S_{\mu ,\nu } \left( z \right) = \frac{{2^{\mu  + 1} }}{{\Gamma \left( {\frac{{\nu  - \mu  + 1}}{2}} \right)\Gamma \left( {\frac{{1 - \mu  - \nu }}{2}} \right)}}\pi e^{ - \frac{\pi }{2}i\mu } K_\nu  \left( {iz} \right) - e^{ - \pi i\mu } S_{\mu ,\nu } \left( {ze^{\pi i} } \right)
\]
for $-\frac{3\pi}{2} < \arg z < -\frac{\pi}{2}$. For the right-hand sides, we can apply the large-$z$ asymptotic expansions of the modified Bessel function and the Lommel function to deduce that
\[
S_{\mu ,\nu } \left( z \right) \sim \frac{{2^{\mu  + 1} }}{{\Gamma \left( {\frac{{\nu  - \mu  + 1}}{2}} \right)\Gamma \left( {\frac{{1 - \mu  - \nu }}{2}} \right)}}\pi e^{\frac{\pi }{2}i\mu } \left( {\frac{\pi }{{ - 2iz}}} \right)^{\frac{1}{2}} e^{iz} \sum\limits_{m = 0}^\infty  {\frac{{a_m \left( \nu  \right)}}{{\left( { - iz} \right)^m }}}  + z^{\mu  - 1} \sum\limits_{n = 0}^\infty  {\left( { - 1} \right)^n \frac{{a_n \left( { - \mu ,\nu } \right)}}{{z^{2n} }}} 
\]
as $z\to \infty$ in the sector $\frac{\pi}{2} < \arg z < \frac{3\pi}{2}$, and
\[
S_{\mu ,\nu } \left( z \right) \sim \frac{{2^{\mu  + 1} }}{{\Gamma \left( {\frac{{\nu  - \mu  + 1}}{2}} \right)\Gamma \left( {\frac{{1 - \mu  - \nu }}{2}} \right)}}\pi e^{ - \frac{\pi }{2}i\mu } \left( {\frac{\pi }{{2iz}}} \right)^{\frac{1}{2}} e^{ - iz} \sum\limits_{m = 0}^\infty  {\frac{{a_m \left( \nu  \right)}}{{\left( {iz} \right)^m }}}  + z^{\mu  - 1} \sum\limits_{n = 0}^\infty  {\left( { - 1} \right)^n \frac{{a_n \left( { - \mu ,\nu } \right)}}{{z^{2n} }}}
\]
as $z \to \infty$ in the sector $-\frac{3\pi}{2} < \arg z < -\frac{\pi}{2}$. Therefore, as the line $\arg z = \frac{\pi}{2}$ is crossed, the additional series
\begin{equation}\label{eq19}
\frac{{2^{\mu  + 1} }}{{\Gamma \left( {\frac{{\nu  - \mu  + 1}}{2}} \right)\Gamma \left( {\frac{{1 - \mu  - \nu }}{2}} \right)}}\pi e^{\frac{\pi }{2}i\mu } \left( {\frac{\pi }{{ - 2iz}}} \right)^{\frac{1}{2}} e^{iz} \sum\limits_{m = 0}^\infty  {\frac{{a_m \left( \nu  \right)}}{{\left( { - iz} \right)^m }}}
\end{equation}
appears in the asymptotic expansion of $S_{\mu ,\nu } \left( z \right)$ beside the original one \eqref{eq15}. Similarly, as we pass through the line $\arg z = -\frac{\pi}{2}$, the series
\begin{equation}\label{eq20}
\frac{{2^{\mu  + 1} }}{{\Gamma \left( {\frac{{\nu  - \mu  + 1}}{2}} \right)\Gamma \left( {\frac{{1 - \mu  - \nu }}{2}} \right)}}\pi e^{ - \frac{\pi }{2}i\mu } \left( {\frac{\pi }{{2iz}}} \right)^{\frac{1}{2}} e^{ - iz} \sum\limits_{m = 0}^\infty  {\frac{{a_m \left( \nu  \right)}}{{\left( {iz} \right)^m }}}
\end{equation}
appears in the asymptotic expansion of $S_{\mu ,\nu } \left( z \right)$ beside the original series \eqref{eq15}. We have encountered a Stokes phenomenon with Stokes lines $\arg z = \pm\frac{\pi}{2}$.

In the important papers \cite{Berry3, Berry2}, Berry provided a new interpretation of the Stokes phenomenon; he found that assuming optimal truncation, the transition between compound asymptotic expansions is of Error function type, thus yielding a smooth, although very rapid, transition as a Stokes line is crossed.

Using the exponentially improved expansion given in Theorem \ref{thm2}, we show that the asymptotic expansion of $S_{\mu ,\nu } \left( z \right)$ exhibits the Berry transition between the two asymptotic series across the Stokes lines $\arg z = \pm\frac{\pi}{2}$. More precisely, we shall find that the first few terms of the series in \eqref{eq19} and \eqref{eq20} ``emerge" in a rapid and smooth way as $\arg z$ passes through $\frac{\pi}{2}$ and $-\frac{\pi}{2}$, respectively.

From Theorem \ref{thm2}, we conclude that if $N \approx \frac{1}{2}\left| z \right|$, then for large $z$, $ \left|\arg z\right| < \pi$, we have
\begin{align*}
S_{\mu ,\nu } \left( z \right) \approx & \; z^{\mu  - 1} \sum\limits_{n = 0}^{N-1}  {\left( { - 1} \right)^n \frac{{a_n \left( { - \mu ,\nu } \right)}}{{z^{2n} }}} \\ & + \frac{{2^{\mu  + 1} }}{{\Gamma \left( {\frac{{\nu  - \mu  + 1}}{2}} \right)\Gamma \left( {\frac{{1 - \mu  - \nu }}{2}} \right)}}\pi e^{\frac{\pi }{2}i\mu } \left( {\frac{\pi }{{ - 2iz}}} \right)^{\frac{1}{2}} e^{iz} \sum\limits_{m = 0} {\frac{{a_m \left( \nu  \right)}}{{\left( { - iz} \right)^m }}\widehat T_{2N - m - \mu  + \frac{1}{2}} \left( {iz} \right)} 
\\ & + \frac{{2^{\mu  + 1} }}{{\Gamma \left( {\frac{{\nu  - \mu  + 1}}{2}} \right)\Gamma \left( {\frac{{1 - \mu  - \nu }}{2}} \right)}}\pi e^{ - \frac{\pi }{2}i\mu } \left( {\frac{\pi }{{2iz}}} \right)^{\frac{1}{2}} e^{ - iz} \sum\limits_{m = 0} {\frac{{a_m \left( \nu  \right)}}{{\left( {iz} \right)^m }}e^{2\pi i\mu } \widehat T_{2N - m - \mu  + \frac{1}{2}} \left( { - iz} \right)} 
\end{align*}
where $\sum\nolimits_{m = 0}$ means that the sum is restricted to the first few terms of the series.

In the upper half-plane the terms involving $\widehat T_{2N - m -\mu + \frac{1}{2}} \left( { - i z} \right)$ are exponentially small, the dominant contribution comes from the terms involving $\widehat T_{2N - m -\mu + \frac{1}{2}} \left( {iz} \right)$. Under the above assumption on $N$, from \eqref{eq12} and \eqref{eq14}, the Terminant functions have the asymptotic behaviour
\[
\widehat T_{2N - m -\mu  + \frac{1}{2}} \left( {iz} \right) \sim \frac{1}{2} + \frac{1}{2}\mathop{\text{erf}}\left( {\left( {\theta  - \frac{\pi }{2}} \right)\sqrt {\frac{1}{2}\left| z \right|} } \right)
\]
provided that $\arg z = \theta$ is close to $\frac{\pi}{2}$, $z$ is large, $\mu$ and $m$ are small in comparison with $N$. Therefore, when $\theta  < \frac{\pi}{2}$, the Terminant functions are exponentially small; for $\theta  = \frac{\pi }{2}$, they are asymptotically $\frac{1}{2}$ up to an exponentially small error; and when $\theta  >  \frac{\pi}{2}$, the Terminant functions are asymptotic to $1$ with an exponentially small error. Thus, the transition across the Stokes line $\arg z = \frac{\pi}{2}$ is effected rapidly and smoothly. Similarly, in the lower half-plane, the dominant contribution is controlled by the terms involving $\widehat T_{2N - m -\mu + \frac{1}{2}} \left( { - iz} \right)$. From \eqref{eq13} and \eqref{eq14}, we have
\[
e^{2\pi i\mu } \widehat T_{2N - m - \mu  + \frac{1}{2}} \left( { - iz} \right) \sim \frac{1}{2} - \frac{1}{2}\mathop{\text{erf}}\left( {\left( {\theta  + \frac{\pi }{2}} \right)\sqrt {\frac{1}{2}\left| z \right|} } \right)
\]
under the assumptions that $\arg z = \theta$ is close to $-\frac{\pi}{2}$, $z$ is large, $\mu$ and $m$ are small in comparison with $N \approx \frac{1}{2}\left| z  \right|$. Thus, when $\theta  >  - \frac{\pi}{2}$, the normalised Terminant functions are exponentially small; for $\theta  =  -\frac{\pi}{2}$, they are asymptotic to $\frac{1}{2}$ with an exponentially small error; and when $\theta < - \frac{\pi}{2}$, the normalised Terminant functions are asymptotically $1$ up to an exponentially small error. Therefore, the transition through the Stokes line $\arg z = -\frac{\pi}{2}$ is carried out rapidly and smoothly.

\subsection{Proof of Theorem \ref{thm4}}\label{subsection43}

Let us denote $z = re^{i\theta }$ with $\left|\theta\right|<\frac{\pi}{2}$ and $r>0$. The remainder $R_N \left( {z,\mu ,\nu } \right)$ may be written
\begin{equation}\label{eq37}
R_N \left( {z,\mu ,\nu } \right) = \left( { - 1} \right)^N \frac{{2^{\mu  + 1} \left( {e^{i\theta } } \right)^{\mu  - 2N - 1} }}{{\Gamma \left( {\frac{{\nu  - \mu  + 1}}{2}} \right)\Gamma \left( {\frac{{1 - \mu  - \nu }}{2}} \right)}}\int_0^{ + \infty } {\frac{{\tau ^{2N - \mu } K_\nu  \left( {r\tau } \right)}}{{1 + \tau ^2 e^{ - 2i\theta } }}d\tau } 
\end{equation}
provided that $\Re \left( \mu  \right) + \left| {\Re \left( \nu  \right)} \right| \le 2N + 1$. The asymptotic expansion of the Lommel function is obtained by expanding the slowly varying part $1/\left(1 + \tau ^2 e^{ - 2i\theta }\right)$ around $\tau=0$ and integrating the resulting series term-by-term. However, when $N$ becomes comparable with $r$, the factor $\tau ^{2N - \mu } K_\nu  \left( {r\tau } \right)$ is dominated by its behaviour near $\tau  = 2N/r$ and not near $\tau=0$. Therefore, it is reasonable to expect that if the asymptotic expansion \eqref{eq1} is truncated when $n\approx r/2$, and the factor $1/\left(1 + \tau ^2 e^{ - 2i\theta }\right)$ in \eqref{eq37} as a function of $\tau^2$ is expanded around $1$, the resulting series provides much better approximations than the original asymptotic expansion \eqref{eq1}.

Let $M$ be an arbitrary non-negative integer and set $\alpha = 1/\left(1 + e^{2i\theta }\right)$, then
\[
\frac{1}{{1 + \tau ^2 e^{ - 2i\theta } }} = \alpha e^{2i\theta } \sum\limits_{m = 0}^{M - 1} {\alpha ^m \left( {1 - \tau ^2 } \right)^m }  + \frac{{\alpha ^M \left( {1 - \tau ^2 } \right)^M }}{{1 + \tau ^2 e^{ - 2i\theta } }}.
\]
Substitution into \eqref{eq37} yields
\begin{gather}\label{eq38}
\begin{split}
R_N \left( {z,\mu ,\nu } \right) = \; & \left( { - 1} \right)^N \frac{{2^{\mu  + 1} \left( {e^{i\theta } } \right)^{\mu  - 2N - 1} }}{{\Gamma \left( {\frac{{\nu  - \mu  + 1}}{2}} \right)\Gamma \left( {\frac{{1 - \mu  - \nu }}{2}} \right)}}\alpha e^{2i\theta } \sum\limits_{m = 0}^{M - 1} {\alpha ^m \int_0^{ + \infty } {K_\nu  \left( {r\tau } \right)\tau ^{2N - \mu } \left( {1 - \tau ^2 } \right)^m d\tau } } \\ & + \widehat R_{N,M} \left( {z,\mu ,\nu } \right),
\end{split}
\end{gather}
where
\begin{equation}\label{eq39}
\widehat R_{N,M} \left( {z,\mu ,\nu } \right) = \left( { - 1} \right)^N \frac{{2^{\mu  + 1} \left( {e^{i\theta } } \right)^{\mu  - 2N - 1}}}{{\Gamma \left( {\frac{{\nu  - \mu  + 1}}{2}} \right)\Gamma \left( {\frac{{1 - \mu  - \nu }}{2}} \right)}}\alpha ^M \int_0^{ + \infty } {\frac{{K_\nu  \left( {r\tau } \right)\tau ^{2N - \mu } \left( {1 - \tau ^2 } \right)^M }}{{1 + \tau ^2 e^{ - 2i\theta } }}d\tau } .
\end{equation}
By expanding the factor $\left(1-\tau^2\right)^m$ in \eqref{eq38}, we derive
\begin{equation}\label{eq41}
R_{N} \left( {z,\mu ,\nu } \right) =z^{\mu - 1} \sum\limits_{m = 0}^{M - 1} {v_{N,m} \left( {z,\mu ,\nu } \right)}  + \widehat R_{N,M} \left( {z,\mu ,\nu } \right),
\end{equation}
where
\begin{align*}
v_{N,m} \left( {z,\mu ,\nu } \right) & = \frac{\left( { - 1} \right)^N}{z^{2N}} \frac{{e^{ - 2i\theta m} }}{{\left( {1 + e^{ - 2i\theta } } \right)^{m + 1} }}\sum\limits_{k = 0}^m {\binom{m}{k}\left( { - 1} \right)^k \frac{{a_{N + k} \left( { - \mu ,\nu } \right)}}{{\left| z \right|^{2k} }}}
\\ & = \frac{1}{{\left( {1 + e^{ - 2i\theta } } \right)^{m + 1} }}\sum\limits_{k = 0}^m {\binom{m}{k}\left( { - 1} \right)^{N + k} \frac{{a_{N + k} \left( { - \mu ,\nu } \right)}}{{z^{2N + 2k} }}e^{ - 2i\theta \left( {m - k} \right)} } .
\end{align*}
Naturally, all these manipulations are valid without the assumption $N \approx r/2$, however, as we shall see later, the expansion \eqref{eq41} is exponentially accurate if $N$ is close to $r/2$.

If $M$ is even, we can derive a simple bound for the error term $\widehat R_{N,M} \left( {z,\mu ,\nu } \right)$. In order to avoid long and complicated expressions, we introduce the notation
\[
\widehat\ell \left( \theta  \right) = \begin{cases} \left|\csc \theta \right| & \text{ if } \; \frac{\pi}{4} < \left|\theta\right| <\frac{\pi}{2} \\ 2\left| {\cos\theta } \right| & \text{ if } \; \left|\theta\right| \leq \frac{\pi}{4} . \end{cases} 
\]
Simple estimation of \eqref{eq39}, using $\ell \left( \theta  \right)$ introduced in \eqref{eq40}, gives us the error bound
\begin{align*}
\left| {\widehat R_{N,M} \left( {z,\mu ,\nu } \right)} \right| & \le \frac{{2^{\Re\left(\mu\right)  + 1} }}{{\left| {\Gamma \left( {\frac{{\nu  - \mu  + 1}}{2}} \right)\Gamma \left( {\frac{{1 - \mu  - \nu }}{2}} \right)} \right|}}\left| \alpha  \right|^M \ell \left( \theta  \right)\int_0^{ + \infty } {K_{\Re \left( \nu  \right)} \left( {r\tau } \right)\tau ^{2N - \Re \left( \mu  \right)} \left( {1 - \tau ^2 } \right)^M d\tau } 
\\ & = \left| {\frac{{\Gamma \left( {\frac{{\Re \left( \nu  \right) - \Re \left( \mu  \right) + 1}}{2}} \right)\Gamma \left( {\frac{{1 - \Re \left( \mu  \right) - \Re \left( \nu  \right)}}{2}} \right)}}{{\Gamma \left( {\frac{{\nu  - \mu  + 1}}{2}} \right)\Gamma \left( {\frac{{1 - \mu  - \nu }}{2}} \right)}}} \right|\left| z \right|^{\Re \left( \mu  \right) - 1} \left| {v_{N,M} \left( {z,\Re \left( \mu  \right),\Re \left( \nu  \right)} \right)} \right|\widehat\ell \left( \theta  \right).
\end{align*}
If $\Re \left( \nu  \right) + \Re \left( \mu  \right) \neq \nu + \mu$ and $\Re \left( \nu  \right) - \Re \left( \mu  \right) \neq \nu - \mu$, and neither number on the right-hand side is a negative odd integer, then we can employ the inequality \eqref{eq36} to simplify this estimate to
\[
\left| {\widehat R_{N,M} \left( {z,\mu ,\nu } \right)} \right| \le \left| {\frac{{\cos \left( {\pi \mu } \right) + \cos \left( {\pi \nu } \right)}}{{\cos \left( {\pi \Re \left( \mu  \right)} \right) + \cos \left( {\pi \Re \left( \nu  \right)} \right)}}} \right|\left| z \right|^{\Re \left( \mu  \right) - 1} \left| {v_{N,M} \left( {z,\Re \left( \mu  \right),\Re \left( \nu  \right)} \right)} \right|\widehat\ell \left( \theta  \right).
\]
If $\Re \left( \nu  \right) + \Re \left( \mu  \right) \neq \nu + \mu$ or $\Re \left( \nu  \right) - \Re \left( \mu  \right) \neq \nu - \mu$ and at least one of them is an positive odd integer, then the limiting value has to be taken in this estimate. The existence of the limit follows from the definition of the coefficients $a_{N+k} \left( { - \Re \left( \mu  \right),\Re \left( \nu  \right)} \right)$ and the fact that $\left| {\Re \left( \nu  \right)} \right| + \Re \left( \mu  \right) < 2N + 1$. In particular, if  both $\mu$ and $\nu$ are real, the absolute value of the remainder term $\widehat R_{N,M} \left( {z,\mu ,\nu } \right)$ is bounded by at most twice the absolute value of the first omitted term of the series \eqref{eq41}, provided that $\left| \nu  \right| + \mu  < 2N + 1$ and that $M$ is even.

Now assume that $N=\frac{1}{2}\left|z\right| +\rho$ and $\left|\arg z\right|\leq \frac{\pi}{2}-\delta <\frac{\pi}{2}$ with some fixed $\rho$ and $0<\delta \leq \frac{\pi}{2}$. As before, $M$ is an arbitrary fixed non-negative integer, not necessarily even. Simple estimation of \eqref{eq39} yields the bound
\[
\left| {\widehat R_{N,M} \left( {z,\mu ,\nu } \right)} \right| \le \frac{{2^{\Re \left( \mu  \right) + 1} }}{{\left| {\Gamma \left( {\frac{{\nu  - \mu  + 1}}{2}} \right)\Gamma \left( {\frac{{1 - \mu  - \nu }}{2}} \right)} \right|}}\frac{{\ell \left( \theta  \right)}}{{\left( {2\left| {\sin \delta } \right|} \right)^M }}\int_0^{ + \infty } {K_{\Re \left( \nu  \right)} \left( {r\tau } \right)\tau ^{2N - \Re \left( \mu  \right)} \left| {1 - \tau ^2 } \right|^M d\tau } .
\]
Noting that
\[
\ell \left( \theta  \right) \leq \begin{cases} \left|\csc \left(2\delta\right) \right| & \text{ if } \; \frac{\pi}{4} < \left|\theta\right| \leq \frac{\pi}{2}-\delta \\ 1 & \text{ if } \; \left|\theta\right| \leq \frac{\pi}{4} , \end{cases} 
\]
we obtain the simple order estimate
\begin{equation}\label{eq51}
\widehat R_{N,M} \left( {z,\mu ,\nu } \right) = \mathcal{O}_{M,\mu ,\nu ,\delta } \left( {\int_0^{ + \infty } {K_{\Re \left( \nu  \right)} \left( {r\tau } \right)\tau ^{2N - \Re \left( \mu  \right)} \left| {1 - \tau ^2 } \right|^M d\tau } } \right).
\end{equation}
It remains to consider the asymptotic behaviour of the integral under the $\mathcal{O}$-symbol for large $\left|z\right|$. Since $N=\frac{1}{2}r +\rho$, we can write this integral as
\begin{multline*}
\int_0^{1/r} {K_{\Re \left( \nu  \right)} \left( {r\tau } \right)\tau ^{r + 2\rho  - \Re \left( \mu  \right)} \left( {1 - \tau ^2 } \right)^M d\tau }  + \int_{1/r}^1 {K_{\Re \left( \nu  \right)} \left( {r\tau } \right)\tau ^{r + 2\rho  - \Re \left( \mu  \right)} \left( {1 - \tau ^2 } \right)^M d\tau } \\ + \int_1^{ + \infty } {K_{\Re \left( \nu  \right)} \left( {r\tau } \right)\tau ^{r + 2\rho  - \Re \left( \mu  \right)} \left( {\tau ^2  - 1} \right)^M d\tau } .
\end{multline*}
The asymptotic behaviour of the modified Bessel function near $0$ shows that the first integral is
\begin{equation}\label{eq50}
\mathcal{O}_{\mu ,\nu } \left( {\frac{1}{{\left| z \right|^{\left| z \right| + 2\rho  - \Re \left( \mu  \right) + 2} }}} \right),
\end{equation}
for large $\left|z\right|$. To estimate the second two integrals, we use the fact that for $x\geq 1$
\[
K_\nu  \left( x \right) = \mathcal{O}_\nu  \left( {\frac{{e^{ - x} }}{{\sqrt {2\pi x} }}} \right).
\]
This follows from the known error bounds for the large-$x$ asymptotics of the modified Bessel function (see, e.g., \cite[10.40.iii]{NIST}). Substituting this expression into the second two integrals and using Laplace's method \cite[pp. 80--82]{Olver}, their contribution is found to be
\[
\mathcal{O}_{M,\mu ,\nu ,\rho } \left( {\frac{{e^{ - \left| z \right|} }}{{\left| z \right|^{\frac{M}{2} + 1} }}} \right),
\]
for large $\left|z\right|$. Clearly, this contribution dominates over \eqref{eq50} as $\left|z\right|\to +\infty$, and together with \eqref{eq51} implies \eqref{eq52}.

\section{Asymptotic expansions for the Struve function}\label{section5}

The formulas \eqref{eq44}, \eqref{eq45} and \eqref{eq5} imply
\begin{equation}\label{eq46}
\mathbf{M}_\nu  \left( z \right) = \frac{1}{\pi }\sum\limits_{n = 0}^{N - 1} {\left( { - 1} \right)^{n + 1} \frac{{\Gamma \left( {n + \frac{1}{2}} \right)\left( {\frac{1}{2}z} \right)^{\nu  - 2n - 1} }}{{\Gamma \left( {\nu  - n + \frac{1}{2}} \right)}}}  + \left( { - 1} \right)^{N + 1} \frac{2}{\pi }z^{\nu  - 2N - 1} \int_0^{ + \infty } {\frac{{t^{2N - \nu } J_\nu  \left( t \right)}}{{1 + \left( {t/z} \right)^2 }}dt} ,
\end{equation}
provided that $\left|\arg z\right|<\frac{\pi}{2}$ and $2N - \frac{3}{2} < \Re \left( \nu  \right)$. If we neglect the term involving the integral and extend the sum to $N=\infty$, we recover the known large-$z$ asymptotic series of the Struve function $\mathbf{M}_\nu  \left( z \right)$. In what follows, we assume that $\nu>0$. Our aim is to show that the large-$z$ asymptotic expansion serves as a uniform large order asymptotic series as well. This interesting property of the large-$z$ asymptotic series of $\mathbf{M}_\nu  \left( z \right)$ was already observed by Dingle \cite[pp. 389--391]{Dingle}, though he did not provide a rigorous proof. Let $\lambda$ be an arbitrary positive real number and take $z=\lambda \nu$ in \eqref{eq46} to find
\begin{multline*}
\mathbf{M}_\nu  \left( {\lambda \nu } \right) =  - \frac{{\left( {\frac{1}{2}\lambda \nu } \right)^{\nu  - 1} }}{{\sqrt \pi  \Gamma \left( {\nu  + \frac{1}{2}} \right)}}\left( {\sum\limits_{n = 0}^{N - 1} {\left( { - 1} \right)^n \frac{\left(2n\right)!}{n!  \lambda ^{2n} }\varphi _n \left( \nu  \right)} }\right. \\ \left.{+ \left( { - 1} \right)^N \frac{{\left( {2N} \right)!}}{{N!\lambda ^{2N} }}\varphi _N \left( \nu  \right)\frac{{2^{\nu  - 2N} \Gamma \left( {\nu  - N + \frac{1}{2}} \right)}}{{\Gamma \left( {N + \frac{1}{2}} \right)}}\int_0^{ + \infty } {\frac{{t^{2N - \nu } J_\nu  \left( t \right)}}{{1 + \left( {t/\lambda \nu } \right)^2 }}dt} } \right),
\end{multline*}
for $\max\left(0,2N - \frac{3}{2}\right) < \nu$, where
\[
\varphi _n \left( \nu  \right) = \frac{{\Gamma \left( {\nu  + \frac{1}{2}} \right)}}{{\nu ^{2n} \Gamma \left( {\nu  - n + \frac{1}{2}} \right)}} \sim \frac{1}{\nu^n}
\]
forms an asymptotic sequence as $\nu \to +\infty$. In Appendix \ref{appendixa}, we show that
\begin{equation}\label{eq43}
\frac{{2^{\nu  - 2N} \Gamma \left( {\nu  - N + \frac{1}{2}} \right)}}{{\Gamma \left( {N + \frac{1}{2}} \right)}}\int_0^{ + \infty } {\frac{{t^{2N - \nu } J_\nu  \left( t \right)}}{{1 + \left( {t/\lambda \nu } \right)^2 }}dt}  = \mathcal{O}\left( 1 \right)
\end{equation}
as $\nu \to +\infty$, uniformly with respect to $\lambda$. Therefore, for any fixed $\lambda>0$, we have the generalised asymptotic expansion
\begin{equation}\label{eq47}
\mathbf{M}_\nu  \left( {\lambda \nu } \right) \sim  - \frac{{\left( {\frac{1}{2}\lambda \nu } \right)^{\nu  - 1} }}{{\sqrt \pi  \Gamma \left( {\nu  + \frac{1}{2}} \right)}}\sum\limits_{n = 0}^\infty  {\left( { - 1} \right)^n \frac{\left(2n\right)!}{n!  \lambda ^{2n} }\varphi _n \left( \nu  \right)} 
\end{equation}
as $\nu \to +\infty$, uniformly with respect to $\lambda$. Interestingly, this asymptotic series provides an explicit formula for the coefficients of the standard large order asymptotic expansion of the Struve function \cite[11.6.E7]{NIST}
\begin{equation}\label{eq49}
\mathbf{M}_\nu  \left( {\lambda \nu } \right) \sim  - \frac{{\left( {\frac{1}{2}\lambda \nu } \right)^{\nu  - 1} }}{{\sqrt \pi  \Gamma \left( {\nu  + \frac{1}{2}} \right)}} \sum\limits_{n = 0}^\infty  {\frac{{n!c_n \left( {i\lambda } \right)}}{{\nu ^n }}} .
\end{equation}
Here, $c_n\left(\lambda\right)$ is a polynomial in $\lambda^{-2}$ of degree $n$. From the asymptotic expansion of the ratio of two Gamma functions \cite[5.11.E13]{NIST}, one finds the exact expression
\begin{equation}\label{eq48}
\varphi _n \left( \nu  \right) = \frac{1}{{\nu ^n }}\sum\limits_{k = 0}^n {\binom{n}{k}\frac{{B_k^{\left( {n + 1} \right)} \left( {\frac{1}{2}} \right)}}{{\nu ^k }}}, 
\end{equation}
where $B_k^{\left( \kappa  \right)} \left(\lambda\right)$ denotes the generalised Bernoulli polynomials, which are defined by the exponential generating function
\[
\left( \frac{z}{e^z  - 1} \right)^\kappa  e^{\lambda z}  = \sum\limits_{k = 0}^\infty  {B_k^{\left( \kappa  \right)} \left(\lambda\right)\frac{z^k}{k!}} \; \text{ for } \; \left|z\right| < 2\pi.
\]
For basic properties of these polynomials, see Milne-Thomson \cite{Milne-Thomson} or N\"{o}rlund \cite{Norlund}. Substituting \eqref{eq48} into \eqref{eq47} and expanding the series in inverse powers of $\nu$, one finds the simple explicit formula
\[
c_n \left( \lambda  \right) = \sum\limits_{k = \left\lceil {n/2} \right\rceil }^n {\binom{2k}{n}\frac{{B_{n - k}^{\left( {k + 1} \right)} \left( {\frac{1}{2}} \right)}}{{\left( {n - k} \right)!}}\frac{1}{\lambda ^{2k}}} .
\]
The first few coefficients are found to be
\[
c_0 \left( \lambda  \right) = 1,\; c_1 \left( \lambda  \right) = 2\lambda ^{ - 2} ,\; c_2 \left( \lambda  \right) = 6\lambda ^{ - 4}  - \frac{1}{2}\lambda ^{ - 2} ,\; c_3 \left( \lambda  \right) = 20\lambda ^{ - 6}  - 4\lambda ^{ - 4} .
\]
These are in agreement with those given in \cite[11.6.E8]{NIST}.

We remark that \eqref{eq49} is valid for $\left|\arg \nu \right| \leq \frac{\pi}{2}-\delta<\frac{\pi}{2}$ with a fixed $0<\delta \leq \frac{\pi}{2}$, and by rearranging its terms, we find that \eqref{eq47} remains valid in the wider range $\left|\arg \nu \right| \leq \frac{\pi}{2}-\delta<\frac{\pi}{2}$.

\appendix

\section{Converging factor}\label{appendixa}

In the sector $\left|\arg w\right| \leq \pi -\delta < \pi$, $0 < \delta \leq \pi$, the Terminant function $\widehat T_p \left( w \right)$ has the asymptotic expansion
\begin{equation}\label{eq22}
\widehat T_p \left( w \right) \sim  - \frac{{ie^{\left( {\pi  -\varphi } \right)pi} }}{{1 + e^{ - i\varphi } }}\frac{{e^{ - w - \left| w \right|} }}{{\sqrt {2\pi \left| w \right|} }}\sum\limits_{k = 0}^\infty  {\frac{{c_k \left( {\varphi ,p - \left| w \right|} \right)}}{{\left| w \right|^k }}} 
\end{equation}
as $w\to +\infty$ provided that $\left|p\right| \sim \left|w\right|$. The coefficients $c_k \left( {\varphi ,p - \left| w \right|} \right)$ are polynomials in $\left(1 + e^{i\varphi }\right)^{-1}$ and $p - \left| w \right|$, respectively \cite[p. 261]{Paris3}. Define the converging factor $\mathscr{C}_N \left( {z,\mu ,\nu } \right)$ via the expression
\[
S_{\mu,\nu}\left(z\right) = z^{\mu  - 1} \sum\limits_{n = 0}^{N - 1} {\left( { - 1} \right)^n \frac{{a_n \left( { - \mu ,\nu } \right)}}{{z^{2n} }}}  + z^{\mu  - 1} \left( { - 1} \right)^N \frac{{a_N \left( { - \mu ,\nu } \right)}}{{z^{2N} }}\mathscr{C}_N \left( {z,\mu ,\nu } \right).
\]
Suppose that $z=\left( {2N + \zeta } \right)e^{i\theta }$ with a fixed real $\zeta$ and $\left|\theta\right|\leq \frac{\pi}{2}-\delta<\frac{\pi}{2}$, $0<\delta \leq \frac{\pi}{2}$. Using the relation $\mathscr{C}_N \left( {z,\mu ,\nu } \right) = \left( { - 1} \right)^N z^{2N - \mu  + 1} R_N \left( {z,\mu ,\nu } \right)/a_N \left( { - \mu ,\nu } \right)$, Theorem \ref{thm2} and the asymptotic expansion \eqref{eq22}, it follows that the converging factor has the asymptotic series
\begin{equation}\label{eq23}
\mathscr{C}_N \left( {z,\mu ,\nu } \right) \sim \frac{{\left| z \right|^{2N - \mu } e^{ - \left| z \right|} \pi }}{{2^{2N - \mu } \Gamma \left( {N + \frac{{\nu  - \mu  + 1}}{2}} \right)\Gamma \left( {N+\frac{{1 - \mu  - \nu }}{2}} \right)}}\sum\limits_{n = 0}^\infty  {\frac{{g_n^ +  \left( {\theta,\zeta} \right) + g_n^ -  \left( {\theta,\zeta} \right)}}{{\left| z \right|^n }}}
\end{equation}
as $z\to \infty$ (or as $N \to +\infty$), with
\[
g_n^ \pm  \left( {\theta,\zeta} \right) = \frac{1}{{1 \mp ie^{ - i\theta } }}\sum\limits_{k = 0}^n {a_k \left( \nu  \right)c_{n - k} \left( {\theta  \pm \frac{\pi }{2},-\zeta - k - \mu  + \frac{1}{2}} \right)}. 
\]
The coefficients $g_n^ \pm  \left( {\theta,\zeta} \right)$ are polynomials in $\left( {1 \pm ie^{i\theta } } \right)^{ - 1}$ and $\zeta$ with coefficients involving $\mu$ and $\nu$. Their sum is a polynomial in $\alpha=\left( {1 + e^{2i\theta } } \right)^{ - 1}$ and $\zeta$. It is well known that for large $N$ and fixed complex $h$, the reciprocal Gamma function possesses the asymptotic expansion
\[
\frac{1}{{\Gamma \left( {N + h} \right)}} \sim \frac{1}{{\sqrt {2\pi } }}N^{ - N - h + \frac{1}{2}} e^N \left( {1 - \frac{{6h^2  - 6h + 1}}{{12N}} +  \cdots } \right) .
\]
Employing this series in \eqref{eq23} for each Gamma function, it is seen that the expansion \eqref{eq23} can be rearranged as an asymptotic series in descending powers of $N$, therefore
\begin{equation}\label{eq24}
\mathscr{C}_N \left( {\left( {2N + \zeta } \right)e^{i\theta } ,\mu ,\nu } \right) \sim \sum\limits_{n = 0}^\infty  {\frac{{\gamma _n \left( {\alpha, \zeta} \right)}}{{N^n }}} ,
\end{equation}
as $N \to +\infty$, where the coefficients $\gamma _n \left( {\alpha, \zeta} \right)$ are polynomials in $\alpha$ and $\zeta$. Note that $\alpha$ becomes unbounded as $\theta \to \pm \frac{\pi}{2}$, whence the sector of validity $\left|\theta\right|\leq \frac{\pi}{2}-\delta<\frac{\pi}{2}$ of this expansion is maximal.

We shall now derive recurrence relations for the polynomial coefficients $\gamma _n \left( {\alpha, \zeta} \right)$ in \eqref{eq24}. From the differential equation satisfied by the Lommel function (see, e.g., \cite[11.9.E1]{NIST}), we obtain that
\begin{equation}\label{eq26}
z^2 \mathscr{C}''_N \left( {z,\mu ,\nu } \right) + z\left( {2\left( {\mu  - 2N - 1} \right) + 1} \right)\mathscr{C}'_N \left( {z,\mu ,\nu } \right) + \left( {\left( {\mu  - 2N - 1} \right)^2  + z^2  - \nu ^2 } \right)\mathscr{C}_N \left( {z,\mu ,\nu } \right) = z^2 ,
\end{equation}
where the differentiation is taken with respect to $z$. Differentiating the series \eqref{eq24} with respect to $\zeta$, we find
\[
e^{i\theta } \mathscr{C}'_N \left( {\left( {2N + \zeta } \right)e^{i\theta } ,\mu ,\nu } \right) \sim \sum\limits_{n = 0}^\infty  {\frac{{\gamma '_n \left( {\alpha ,\zeta } \right)}}{{N^n }}}
\]
and
\[
e^{2i\theta } \mathscr{C}''_N \left( {\left( {2N + \zeta } \right)e^{i\theta } ,\mu ,\nu } \right) \sim \sum\limits_{n = 0}^\infty  {\frac{{\gamma ''_n \left( {\alpha ,\zeta } \right)}}{{N^n }}} ,
\]
as $N\to +\infty$. Substituting into the differential equation \eqref{eq26} and equating the coefficients of powers of $N$, we deduce that
\[
\alpha \gamma ''_0 \left( {\alpha ,\zeta } \right) - 2\alpha \gamma '_0 \left( {\alpha ,\zeta } \right) + \gamma _0 \left( {\alpha ,\zeta } \right) = 1 - \alpha, 
\]
\begin{align*}
& 2\alpha \gamma ''_1 \left( {\alpha ,\zeta } \right) - 4\alpha \gamma '_1 \left( {\alpha ,\zeta } \right) + 2\gamma _1 \left( {\alpha ,\zeta } \right) + 2\alpha \zeta \gamma ''_0 \left( {\alpha ,\zeta } \right) - \alpha \left( {2\zeta  - 2\mu  + 1} \right)\gamma '_0 \left( {\alpha ,\zeta } \right) \\ &+ 2\left( {\zeta \left( {1 - \alpha } \right) + \alpha \left( {1 - \mu } \right)} \right)\gamma _0 \left( {\alpha ,\zeta } \right) = 2\zeta \left( {1 - \alpha } \right),
\end{align*}
\begin{align*}
& 4\alpha \gamma ''_2 \left( {\alpha ,\zeta } \right) - 8\alpha \gamma '_2 \left( {\alpha ,\zeta } \right) + 4\gamma _2 \left( {\alpha ,\zeta } \right) + 4\alpha \zeta \gamma ''_1 \left( {\alpha ,\zeta } \right) - 2\alpha \left( {2\zeta  - 2\mu  + 1} \right)\gamma '_1 \left( {\alpha ,\zeta } \right) \\ &+ 4\left( {\zeta \left( {1 - \alpha } \right) + \alpha \left( {1 - \mu } \right)} \right)\gamma _1 \left( {\alpha ,\zeta } \right) + \alpha \zeta ^2 \gamma ''_0 \left( {\alpha ,\zeta } \right) + \alpha \zeta \left( {2\mu  - 1} \right)\gamma '_0 \left( {\alpha ,\zeta } \right) \\ &+ \left( {\zeta ^2 \left( {1 - \alpha } \right) + \alpha \left( {\left( {\mu  - 1} \right)^2  - \nu ^2 } \right)} \right)\gamma _0 \left( {\alpha ,\zeta } \right) = \zeta ^2 \left( {1 - \alpha } \right)
\end{align*}
and
\begin{gather}\label{eq25}
\begin{split}
& 4\alpha \gamma ''_{n + 2} \left( {\alpha ,\zeta } \right) - 8\alpha \gamma '_{n + 2} \left( {\alpha ,\zeta } \right) + 4\gamma _{n + 2} \left( {\alpha ,\zeta } \right) + 4\alpha \zeta \gamma ''_{n + 1} \left( {\alpha ,\zeta } \right) - 2\alpha \left( {2\zeta  - 2\mu  + 1} \right)\gamma '_{n + 1} \left( {\alpha ,\zeta } \right)\\ & + 4\left( {\zeta \left( {1 - \alpha } \right) + \alpha \left( {1 - \mu } \right)} \right)\gamma _{n + 1} \left( {\alpha ,\zeta } \right) + \alpha \zeta ^2 \gamma ''_n \left( {\alpha ,\zeta } \right) + \alpha \zeta \left( {2\mu  - 1} \right)\gamma '_n \left( {\alpha ,\zeta } \right) \\ & + \left( {\zeta ^2 \left( {1 - \alpha } \right) + \alpha \left( {\left( {\mu  - 1} \right)^2  - \nu ^2 } \right)} \right)\gamma _n \left( {\alpha ,\zeta } \right) =0,
\end{split}
\end{gather}
for $n\geq 1$. The only polynomial solutions of the first three equations are
\[
\gamma _0 \left( {\alpha ,\zeta } \right) = 1 - \alpha ,
\]
\[
\gamma _1 \left( {\alpha ,\zeta } \right) = \alpha \left( {1 - \alpha } \right)\zeta  + \alpha \left( {1 - \alpha } \right)\left( {2\alpha  + \mu  - 1} \right)
\]
and
\begin{align*}
\gamma _2 \left( {\alpha ,\zeta } \right) = \; & \frac{{\alpha \left( {1 - \alpha } \right)\left( {4\alpha  - 3} \right)}}{4}\zeta ^2  + \alpha \left( {1 - \alpha } \right)\left( {6\alpha ^2  + 2\alpha \left( {\mu  - 3} \right) + 1 - \mu } \right)\zeta \\ & + \frac{{\alpha \left( {1 - \alpha } \right)\left( {48\alpha ^3  + 8\alpha ^2 \left( {3\mu  - 8} \right) + 4\alpha \left( {\mu ^2  - 5\mu  + 5} \right) - \left( {\mu  - 1} \right)^2  + \nu ^2 } \right)}}{4}.
\end{align*}
From the fourth equation \eqref{eq25}, it follows by induction that the degree of $\gamma_n \left( {\alpha ,\zeta } \right)$ as a polynomial in $\zeta$ is $n$.

If $z$ is real and positive, the coefficient $\gamma _0 \left( {\alpha ,\zeta } \right)$ is $\frac{1}{2}$. Hence, in this case, if the series \eqref{eq1} is truncated near its smallest term, the remainder is approximately equal to half of the first neglected term.

\section{An auxiliary estimate}\label{appendixb}

We prove the estimate \eqref{eq43}. Let $N$ be a fixed non-negative integer and $\lambda$ be a fixed positive real number. From Stirling's formula, we have
\begin{align*}
\frac{{2^{\nu  - 2N} \Gamma \left( {\nu  - N + \frac{1}{2}} \right)}}{{\Gamma \left( {N + \frac{1}{2}} \right)}}\int_0^{ + \infty } {\frac{{t^{2N - \nu } J_\nu  \left( t \right)}}{{1 + \left( {t/\lambda \nu } \right)^2 }}dt} & = \frac{{2^{\nu  - 2N} \Gamma \left( {\nu  - N + \frac{1}{2}} \right)}}{{\Gamma \left( {N + \frac{1}{2}} \right)}}\nu ^{2N - \nu  + 1} \int_0^{ + \infty } {\frac{{t^{2N - \nu } J_\nu  \left( {\nu t} \right)}}{{1 + \left( {t/\lambda } \right)^2 }}dt}  \\
&  = \mathcal{O}\left( {\left( {\frac{2}{e}} \right)^\nu  \frac{{\nu ^{N + 1} }}{{2^{2N} \Gamma \left( {N + \frac{1}{2}} \right)}}\int_0^{ + \infty } {\frac{{t^{2N - \nu } J_\nu  \left( {\nu t} \right)}}{{1 + \left( {t/\lambda } \right)^2 }}dt} } \right)
\end{align*}
as $\nu \to +\infty$. We split the integral under the $\mathcal{O}$-symbol into three parts and estimate each of them separately. For the interval $0<t<\frac{1}{2}$, we use Watson's inequality \cite[p. 255, expression (9)]{Watson}
\[
\left|J_\nu  \left( {\nu t} \right)\right| \le \frac{{t^\nu  \exp \left( { - \nu \left( {\log \left( {1 + \sqrt {1 - t^2 } } \right) - \sqrt {1 - t^2 } } \right)} \right)}}{{\sqrt {2\pi \nu \sqrt {1 - t^2 } } }},
\]
and Laplace's method \cite[pp. 80--82]{Olver}, to find
\begin{align*}
\left| {\int_0^{\frac{1}{2}} {\frac{{t^{2N - \nu } J_\nu  \left( {\nu t} \right)}}{{1 + \left( {t/\lambda } \right)^2 }}dt} } \right| \le \int_0^{\frac{1}{2}} {t^{2N - \nu } \left|J_\nu  \left( {\nu t} \right)\right|dt} & \le \int_0^{\frac{1}{2}} {t^{2N} \frac{{\exp \left( { - \nu \left( {\log \left( {1 + \sqrt {1 - t^2 } } \right) - \sqrt {1 - t^2 } } \right)} \right)}}{{\sqrt {2\pi \nu \sqrt {1 - t^2 } } }}dt}  \\ &= \mathcal{O}\left( {\left( {\frac{2}{e}} \right)^{ - \nu } \frac{{2^{2N} \Gamma \left( {N + \frac{1}{2}} \right)}}{{\nu ^{N + 1} }}} \right),
\end{align*}
if $\nu \to +\infty$. When $\frac{1}{2} <t<1$, we apply Siegel's bound \cite{Siegel}
\[
\left|J_\nu  \left( {\nu t} \right)\right| \le t^\nu  \exp \left( { - \nu \left( {\log \left( {1 + \sqrt {1 - t^2 } } \right) - \sqrt {1 - t^2 } } \right)} \right),
\]
together with Laplace's method, to deduce that
\begin{align*}
\left| {\int_{\frac{1}{2}}^1 {\frac{{t^{2N - \nu } J_\nu  \left( {\nu t} \right)}}{{1 + \left( {t/\lambda } \right)^2 }}dt} } \right| & \le \int_{\frac{1}{2}}^1 {t^{2N - \nu } \left|J_\nu  \left( {\nu t} \right)\right|dt}  \le \int_{\frac{1}{2}}^1 {t^{2N} \exp \left( { - \nu \left( {\log \left( {1 + \sqrt {1 - t^2 } } \right) - \sqrt {1 - t^2 } } \right)} \right)dt} 
\\ & = \mathcal{O}\left( {\exp \left( { - \nu \left( {\log \left( {1 + \sqrt {\frac{3}{4}} } \right) - \sqrt {\frac{3}{4}} } \right)} \right)\frac{1}{2^{2N}\nu}} \right) = o\left( {\left( {\frac{2}{e}} \right)^{-\nu}  \frac{{2^{2N} \Gamma \left( {N + \frac{1}{2}} \right)}}{\nu ^{N + 1}} } \right)
\end{align*}
as $\nu \to +\infty$. For the remaining case $1 < t<+\infty$, we can use the simple inequality $\left| {J_\nu  \left( {\nu t} \right)} \right| \le 1$ \cite[10.14.E1]{NIST}, to obtain
\[
\int_1^{ + \infty } {\frac{{t^{2N - \nu } J_\nu  \left( {\nu t} \right)}}{{1 + \left( {t/\lambda } \right)^2 }}dt}  = o\left( 1 \right)
\]
if $\nu \to +\infty$, uniformly with respect to $\lambda$. Taking all these estimations together, we deduce that
\begin{multline*}
\frac{{2^{\nu  - 2N} \Gamma \left( {\nu  - N + \frac{1}{2}} \right)}}{{\Gamma \left( {N + \frac{1}{2}} \right)}}\int_0^{ + \infty } {\frac{{t^{2N - \nu } J_\nu  \left( t \right)}}{{1 + \left( {t/\lambda \nu } \right)^2 }}dt} \\ = \mathcal{O}\left( {\left( {\frac{2}{e}} \right)^\nu  \frac{{\nu ^{N + 1} }}{{2^{2N} \Gamma \left( {N + \frac{1}{2}} \right)}}\left( {\mathcal{O}\left( {\left( {\frac{2}{e}} \right)^{ - \nu } \frac{{2^{2N} \Gamma \left( {N + \frac{1}{2}} \right)}}{{\nu ^{N + 1} }}} \right) + o\left( {\left( {\frac{2}{e}} \right)^{ - \nu } \frac{{2^{2N} \Gamma \left( {N + \frac{1}{2}} \right)}}{{\nu ^{N + 1} }}} \right) + o\left( 1 \right)} \right)} \right) = \mathcal{O}\left( 1 \right)
\end{multline*}
when $\nu \to +\infty$, uniformly with respect to $\lambda$.

\end{document}